\documentclass[12pt]{article}
\usepackage{cmap}
\usepackage[T2A]{fontenc}
\usepackage[utf8]{inputenc}
\usepackage[english]{babel}
\usepackage{graphicx}
\graphicspath{{images/}}
\usepackage{wrapfig}
\usepackage{amsmath}
\usepackage{enumitem}
\usepackage{fixltx2e}

\usepackage[backend=biber,style=numeric-comp,maxcitenames=2,maxbibnames=10,isbn=false,doi=false,url=false,firstinits=true]{biblatex}
\newtheorem{remark}{Remark}

\begin{document}

\begin{center}{\large\bf Scenarios for the Creation of Hyperchaotic Attractors \\ of 3D Maps}
\end{center}
\medskip
\begin{center}{\bf Shykhmamedov A.$^1$, Karatetskaia E.$^1$, \\ Kazakov A.$^{1}$, Stankevich N.$^{1}$}
\end{center}
\medskip
\begin{center}{$^1$ National Research University Higher School of Economics,\\
25/12 Bolshaya Pecherskaya Ulitsa, 603155 Nizhny Novgorod, Russia}
\end{center}

\begin{center}
    e-mails: \\ aykhansh@gmail.com; \\ eyukaratetskaya@gmail.com; \\ kazakovdz@yandex.ru; \\ stankevichnv@mail.ru
\end{center}

\begin{abstract}
We study bifurcation mechanisms for the appearance of hyperchaotic attractors in three-dimensional diffeomorphisms, i.e., such attractors whose orbits have two positive Lyapunov exponents in numerical experiments. In order to possess this property periodic orbits belonging to the attractor should have two-dimensional unstable invariant manifolds. For realization of this possibility, we propose several bifurcation scenarios that include cascades of both supercritical period-doubling bifurcations with saddle periodic orbits and supercritical Neimark-Sacker bifurcations with stable periodic orbits, as well as various combinations of these cascades. In the paper, these scenarios are illustrated by an example of the three-dimensional Mir\'a map.
\end{abstract}

{\bf Keywords.} Hyperchaotic attractor, homoclinic orbit, three-dimensional map.

\newpage

\section{Introduction}

This paper is devoted to the study of bifurcation scenarios leading to the appearance of hyperchaotic attractors in three-dimensional map
\begin{equation}
\begin{cases}
\bar x = y, \\
\bar y = z, \\
\bar z = M_1 + Bx + M_2 z - y^2.
\end{cases}
\label{eq:MiraMap}
\end{equation}
Here $x, y$, and $z$ are phase variables; $M_1, M_2$, and $B$ are parameters. This map is the well-known ``homoclinic map'' which was introduced in \cite{GST93b, GST93c} as a normal form of the first return map near homoclinic tangencies of multidimensional systems, see also \cite{GST96, GST08, GGS10}. Map \eqref{eq:MiraMap} has the constant Jacobian $B$. When $B=0$, it becomes effectively two-dimensional map of the form $\bar y = z, \; \bar z = M_1 + M_2 z - y^2$. It is the well-known two-dimensional endomorphism introduced and studied by C. Mir\'a yet in 60s \cite{Mira65}, see also \cite{Mira_book}. Therefore, we call map \eqref{eq:MiraMap} the \textit{three-dimensional Mir\'a map}. It is worth noting that this map belongs the class of 3D H\'enon maps. In the paper \cite{lomeli1998quadratic} by Lomel{\'\i} and Meiss it was shown that every quadratic 3D map with the constant Jacobian, whose the inverse map is also quadratic, can be reduced to a map of the form $\bar x = y, \; \bar y = z, \;  \bar z = Bx + G(y,z)$, where $G$ is a quadratic polynomial in $y$ and $z$. Therefore, such diffeomorphism is also called ``Lomel{\'\i} map'', see e.g. \cite{mireles2013quadratic}.

Recall, that the term ``hyperchaotic attractor'' was introduced by R\"ossler in \cite{Rossler78} for strange attractors with (at least) two directions of exponential instability. In order to be hyperchaotic, an attractor should contain nontrivial hyperbolic subsets whose unstable manifolds have dimension $\geq 2$. Hyperchaoticity of an attractor is simply verified numerically. If numerical experiments show that orbits in the attractor have (at least) two positive Lyapunov exponents (LE), then the observed attractor is hyperchaotic.

As a consequence of the above, the study of hyperchaotic dynamics in map \eqref{eq:MiraMap} looks very promising, since in a certain region of the parameter space (region SH$(1,2)$ in Fig.~\ref{ris:Fig8}) there is a nontrivial hyperbolic set with the simplest structure -- the three-dimensional Smale horseshoe -- in which all its unstable manifolds are two-dimensional \cite{GonLi10}. Hereinafter, we will prescribe a type $(n,m)$ for hyperbolic periodic orbits with $n$-dimensional stable and $m$-dimensional unstable invariant manifolds (where types $(n,0)$ and $(0,m)$ relate, respectively, to the stable and completely unstable periodic orbits). Then, we can denote the above horseshoe as the Smale horseshoe of type (1,2). Accordingly, in the parameter space, along pathways leading from the region where the map has a simple attractor (stable fixed point) to the region where there is such a horseshoe, one can create hyperchaotic attractors. How this is possible and which bifurcations occur along such pathways are the main questions of the current paper.

At first, let us recall some well-known facts on strange attractors in map \eqref{eq:MiraMap}. When $B=0$, the corresponding two-dimensional Mir\'a endomorphism has, in certain regions of the $(M_1,M_2)$-parameter plane, hyperchaotic attractors, the so-called snap-back repellers \cite{Mar78}, consisting mainly of completely unstable orbits (with two positive Lyapunov exponents). As was shown yet in \cite{BaierKlein90}, see also \cite{richter2002generalized}, for small $B$ these attractors are transformed into hyperchaotic attractors for the three-dimensional diffeomorphism \eqref{eq:MiraMap} that keeps a tendency of exponential expansion along two directions. When $B$ is not too small (and $|B|<1$), map \eqref{eq:MiraMap} can demonstrate other types of hyperchaotic attractors, in particular, the so-called \textit{discrete Shilnikov homoclinic attractors} containing a saddle-focus fixed point with two-dimensional unstable invariant manifold \cite{GGS12, GGKT14, GG16}. These attractors are of special interest since they are often observed in multidimensional systems including those from applications.\footnote{In particular, such hyperchaotic attractors were found in the models of gas bubbles' dynamics \cite{GSKK19}, gene represillators \cite{stankevich2021chaos}, in the modified generator of Anischenko-Astakhov \cite{sataev2021cascade} and in many other models. In all these models a discrete Shilnikov attractor appears on the corresponding three- or multidimensional Poincar\'e maps. Moreover, we note that the appearance of such attractors is one of the most typical route to strange attractors in multidimensional systems as whole.} Among the recent studies of map \eqref{eq:MiraMap} we would like to note the works \cite{zhang2016chaotic, zhao2017bifurcation, hampton2022three} where extensive theoretical and numerical researches of regular and chaotic dynamics are presented. In the current paper, we study bifurcation mechanisms for the appearance of hyperchaotic attractors in this map for both cases of small and not very small values of the Jacobian $B$.

In the paper we propose two types of bifurcation scenarios leading to the appearance of hyperchaotic attractors in one-parameter families of three-dimensional maps. In turn, both these scenarios can be divided into two complementary parts. The first part is \textit{phenomenological}, it outlines the main stages of the formation of the so-called discrete \textit{homoclinic attractors} that, by definition \cite{GGS12, GGKT14}, contain only one saddle fixed point and its unstable invariant manifold. In the case under consideration, the fixed point is either a saddle-focus of type (1,2) or a saddle with two negative unstable multipliers. This phenomenological part includes only a few local and global (homoclinic) bifurcations leading from a simple attractor (stable fixed point) to a homoclinic attractor containing nontrivial hyperbolic subsets with two-dimensional unstable manifolds.

However, this attractor may not be hyperchaotic, especially for initial stages of its formation. For the appearance of a strange attractor characterizing by two positive Lyapunov exponents it is necessary that:
\begin{itemize}
\item the majority of orbits in the attractor should have two-dimensional unstable manifolds.\footnote{It is difficult to define the ``majority'' unambiguously here. However, if we assume a small bounded independent and identically distributed (IID) noise acting on the system, then we can expect a unique stationary measure, and the ``majority'' will mean a set of full measure. Note that we do have the round-off noise in the numerical experiments, but it is not clear if this kind of noise is well-modeled by IID noise.}
\end{itemize}
This problem is the essence of the second, \textit{empirical}, part of the scenarios. This part is much more complicated and variable, because, unlike the first one, it includes infinite sequences of various bifurcations. These sequences contain, in particular, bifurcations of periodic orbits (stable and saddle of type (2,1)) transforming them into saddle periodic orbits of type (1,2), giving a way to get two positive Lyapunov exponents in numerical experiments.

The phenomenological part of the scenarios outlines, in general terms, routes to the onset of homoclinic attractors containing a fixed point $O$ of type (1,2). Together with the unstable manifold, such attractors also contain all homoclinic orbits to the fixed point, i.e., such orbits in which the stable and unstable manifolds of the fixed point intersect. In general case, these intersections are transverse. By the Smale-Shilnikov theorem \cite{Sm67,Sh67} this implies that the attractor must contain nontrivial hyperbolic subsets whose unstable invariant manifolds are two-dimensional. However, this fact does not automatically lead to hyperchaos of the attractor, since, as it happens in most cases, this attractor is not hyperbolic and may contain also saddle periodic orbits of type (2,1) (and, even, invisible for numerics periodic sinks) which together can make a greater contribution to the calculation of the averaged Lyapunov exponents.

In particular, this is true for \textit{discrete Shilnikov attractors} containing a saddle-focus fixed point of type (1,2). Such attractors were introduced and partially studied in the papers \cite{GGS12, GGKT14, GG16}, where a phenomenological scenario of their appearance was also described. The main stages of this scenario are schematically represented in Fig.~\ref{ris:Fig6}. They are as follows:
\begin{enumerate}[label=(\roman*)]
\item supercritical Neimark-Sacker bifurcation with the stable fixed point $O$ after which this point becomes saddle-focus of type (1,2) and a stable invariant curve $L$ appears in its neighborhood (Fig.~\ref{ris:Fig6}a);
\item formation of Shilnikov whirlpool tightening almost all orbits from some absorbing domain (Fig.~\ref{ris:Fig6}c);
\item emergence of a transversal intersection between one-dimensional stable and two-dimensional unstable invariant manifolds of the saddle-focus fixed point $O$ (Fig.~\ref{ris:Fig6}d).
\end{enumerate}
In Section~\ref{sec:ShilAttrScen}, we discuss hyperchaotic properties of discrete Shilnikov attractors and study a class of interior bifurcations that create and support the spiral hyperchaos.

We consider also another type of discrete homoclinic attractors which contain a saddle fixed point with a pair of negative unstable multipliers. We show how such attractors can appear in one-parameter families of three-dimensional maps. The key steps of the corresponding phenomenological scenario are as follows:
\begin{enumerate}[label=(\roman*)]
\item supercritical period-doubling bifurcation with a stable fixed point $O$ after which this point becomes saddle of (2,1)-type (Fig.~\ref{ris:Fig2}b);
\item supercritical period-doubling bifurcation with the saddle fixed point $O$ after which it becomes saddle of the desired (2,1)-type (Fig.~\ref{ris:Fig2}c);
\item emergence of a transversal homoclinic intersection between stable and unstable manifolds of $O$ (Fig.~\ref{ris:Fig3}).
\end{enumerate}
Since the proposed new scenario looks as an extension of the well-known scenario of the creation of the H\'enon attractor, we call the corresponding attractors \textit{hyperchaotic H\'enon-like attractors}\footnote{In three-dimensional maps, another type of ``H\'enon-type'' attractors is also possible. They appear as a result of a sequence of period-doubling bifurcations with a stable invariant curve, and were called as \textit{quasiperiodic H\'enon-like attractors} \cite{BSV02, BSV05, BSV10} by analogy with H\'enon-like attractors.}.

In more detail, both phenomenological scenarios are considered in Section~\ref{sec:Scenarios} where we also supplement them by the description of the second, empirical, part responded for the appearance of two positive Lyapunov exponents in numerical experiments. This part is motivated, in a sense, by the papers \cite{KMP2000, YaK01, YaK01_2} where it was shown (on examples of two-dimensional endomorphisms and coupled flow systems) that chaotic attractors can transform to hyperchaotic ones as a result of the absorption of periodic orbits with two-dimensional unstable invariant manifolds. We suggest two bifurcation mechanisms leading to such transformations. The first one is associated with the appearance of infinite cascades of period-doubling bifurcations with saddle periodic orbits of type (2,1) transforming them into orbits of type (1,2). The second mechanism is related to the direct formation of saddle-focus periodic orbits of type (1,2) via Neimark-Sacker bifurcations with stable periodic orbits. In both cases, the absorption of hyperchaotic saddles (of types (1,2)) by the attractor happens via global (homoclinic or heteroclinic) bifurcations.

In the second part of the paper (Sections~\ref{sec:SmallB}--\ref{sec:LargeB}), we apply the proposed scenarios to study mechanisms for the appearance of hyperchaotic attractors in one-parameter families of the three-dimensional Mir\'a map \eqref{eq:MiraMap}. We consider cases of small and not very small values of the Jacobian $B$. In both cases we start from a stable fixed point $O_+$ which appears via a saddle-node bifurcation together with a saddle fixed point $O_{-}$ of type (2,1) and change parameters towards the region SH(1,2), see Figure~\ref{ris:Fig8}, where the Smale horseshoe of type (1,2) exists \cite{GonLi10}.

The first bifurcation on this pathway is always the supercritical Neimark-Sacker bifurcation after which the point $O_+$ becomes saddle-focus of type (1,2) and a stable invariant curve $L$ is born in its neighborhood. Further, this curve breaks down giving, at first, a rise of a certain chaotic attractor (with only one positive Lyapunov exponent) and, then, this attractor becomes hyperchaotic. The study of accompanying bifurcations is, in fact, the main aim of this part of the paper.

We note that usually, before the destruction, the curve $L$ becomes resonant: a pair of period-$q$ stable and saddle (2,1) orbits appears on it (inside the corresponding Arnold tongue); in this case the resonant curve is the closure of the unstable manifold of the period-$q$ saddle orbit.

In the paper we pay special attention to the cases of the so-called \textit{strong resonances 1:3 and 1:4} which are the most difficult and interesting. Besides, their influence on the organization of corresponding bifurcation diagrams is the most visible (see e.g. Lyapunov diagrams for the map under consideration in Fig.~\ref{ris:Fig10}, where codimension-two points giving rise to the Arnold tongues with the strong resonances 1:3 and 1:4 are denoted by R$_3$ and R$_4$, respectively). In the case of two-dimensional diffeomorphisms, the nongenerated 1:3 resonance leads usually to the global instability \cite{SimoVieiro}, the corresponding invariant curve breaks down without the appearance of stable elements of dynamics. However, this is no longer true for three-dimensional maps with not small values of the Jacobian. In particular, as is shown in \cite{GGS12, GGKT14}, regular and chaotic attractors can appear in map \eqref{eq:MiraMap} after the destruction of the curve $L$ near the 1:3 resonance. In this paper, we show that these attractors can be even hyperchaotic. However, they occupy very thin regions in the parameter space, and can appear only for specific (not small) values of the Jacobian $B$, see the red-colored regions below the point R$_3$ in Fig.~\ref{ris:Fig10}d.

As for the 1:4 resonance, we show that, unlike the 1:3 resonance, it responses for the emergence of the most visible Arnold tongue with stable dynamics. Inside this tongue, the resonant period-4 stable and saddle (2,1) orbits undergo numerous bifurcations. Moreover, a type of these bifurcations essentially depends on values of the Jacobian $B$.

For small values of $B$, these bifurcations include a cascade of period-doubling bifurcations resulting in a $4$-component H\'enon-like attractor containing infinitely many saddle orbits of type (2,1), see Fig.~\ref{ris:Fig14}c. In their turn, these orbits, as well as the resonant period-$4$ saddle orbit, also undergoes cascades of period-doubling bifurcations leading to the formation of hyperchaotic attractor containing infinitely many periodic saddle orbits of type (1,2), see Section~\ref{sec:SmallB} for more detail.

The stable period-4 orbit can also undergo the supercritical Neimark-Sacker bifurcation followed by the creation of a four-component discrete Shilnikov-like attractor. This scenario of the transition to hyperchaos is more typical for the cases of not small values of the parameter $B$ ($B \in [0.3,0.6]$ in Fig.~\ref{ris:Fig17}b). Independently of this, the resonant period-4 saddle orbit (of type (1,2)), as in the case of small $B$, undergoes an infinite cascade of period-doubling bifurcations leading to the creation of nontrivial hyperbolic subset of type (1,2). Then, a new type of hyperchaotic attractor can appear when this subset or its parts are absorbed by the above mentioned four-component Shilnikov attractor. After such absorption, the four-component attractor transforms into a one-component attractor containing both period-4 saddle-focus and saddle orbits of type (1,2), see Section~\ref{sec:LargeB}, for more detail.

In the last part of the paper (Section~\ref{sec:VarShil}), we demonstrate another scenario of the destruction of the curve $L$. We show that for not small values of the Jacobian (e.g. for $B=0.7$) this curve can undergo a quite long sequence of period-doubling bifurcations resulting in the formation of a chaotic attractor with one positive, one near-zero (indistinguishable from zero in numerics) and one negative Lyapunov exponents. We discuss and explain this phenomenon.

\section{Scenarios for the appearance of hyperchaotic attractors}\label{sec:Scenarios}

In this section, we describe scenarios (both phenomenological and empirical parts) that lead to the appearance of hyperchaotic attractors in three-dimensional maps. Let us consider a one-parameter family of three-dimensional maps $$\bar x = F(x, \varepsilon)$$ depending on a parameter $\varepsilon$. In the presented below scenarios, we start with a fixed point $O$ which is asymptotically stable and belongs to some absorbing domain $D_a(O)$, see Fig.~\ref{ris:Fig2}a. Finally, as a result of a series of codimension-one bifurcations, we obtain a homoclinic attractor containing the point $O$ which becomes a saddle or a saddle-focus with two-dimensional unstable invariant manifold. We also describe an empirical part of the scenarios due to which the majority of orbits inside an attractor get two-dimensional unstable invariant manifolds.

\subsection{Hyperchaotic H\'enon-like attractor}

The \textit{hyperchaotic H\'enon-like attractor} is a homoclinic attractor containing a saddle fixed point with a pair of negative unstable multipliers. It can appear as a result of the scenario the beginning part of which coincides with the well-known scenario of the birth of the H\'enon attractor. Therefore, let us first recall some details on the H\'enon and H\'enon-like attractors and scenarios for their appearance.

\subsubsection{Some details about H\'enon-like attractors}

\begin{figure}[h!]
\begin{minipage}[h]{1\linewidth}
\center{\includegraphics[width=1\linewidth]{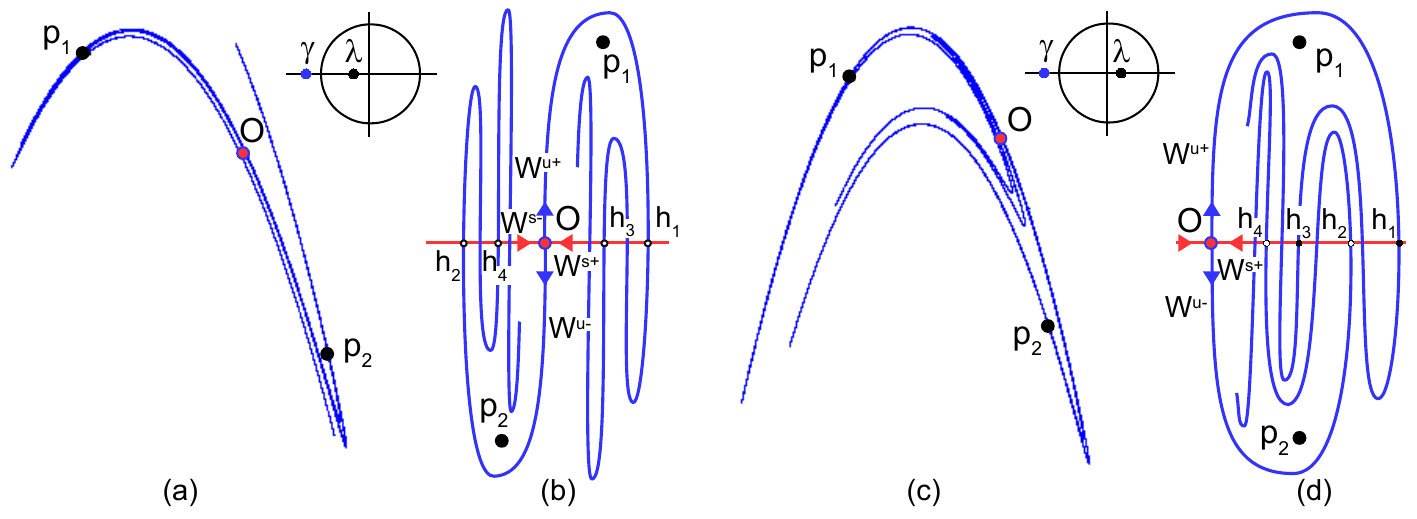}}
\end{minipage}
\caption{\footnotesize (a) Orientable H\'enon attractor in map $H$ for $b = 0.3, M = 2.1$ and (c) its homoclinic structure; (b) nonorientable H\'enon attractor in map $H$ for $b = -0.3, M = 1.4$ and (d) its homoclinic structure. $O$ is a saddle fixed point belonging to the attractor, $W^{s+}, W^{s-}, W^{u+}$, and $W^{u-}$ are pairs of stable and unstable separatrices of $O$, $p_1$ and $p_2$ are components of the period-2 orbit which appears after the period-doubling bifurcation with the point $O$, and $h_1, h_2, h_3, \ldots$ are homoclinic points of some homoclinic to $O$ orbit.}
\label{ris:Fig1}
\end{figure}

The H\'enon attractor is a homoclinic attractor of the two-dimensional H\'enon map $H: \bar x = y, \; \bar y = M - bx - y^2$ \cite{Henon76}. It contains the saddle fixed point $O$ with a negative unstable multiplier and is formed after the Feigenbaum cascade of period-doubling bifurcations \cite{CarLyubMart05, LyubMart11, GambStrTress89} followed by a cascade of heteroclinic ``band-merging'' bifurcations \cite{Simo79}. Let us denote the stable and unstable multipliers of the point $O$ by $\lambda$ and $\gamma$. Then, for the H\'enon attractor, the following conditions always hold: $\gamma < -1$ (the unstable multiplier is negative) and $\sigma = |\gamma\lambda| < 1$ (area-contracting condition). Depending on the sign of the stable multiplier $\lambda$, the H\'enon attractor can be of two types: if $-1 < \lambda < 0$ the attractor is orientable, see the example of its phase portrait in Fig.~\ref{ris:Fig1}a; and if $0< \lambda < 1$ it is nonorientable, see the phase portrait in Fig.~\ref{ris:Fig1}c (exactly this attractor was discovered and studied by M.~H\'enon in \cite{Henon76}).

When dynamics in the map are associated with the existence of the H\'enon attractor\footnote{The question of the existence of the H\'enon attractor is very delicate. For sufficiently small positive values of the Jacobian of the H\'enon map, Benedicks and Carleson proved in \cite{BenCar91} that the set of parameter values corresponding to the chaotic attractor is a Cantor set with the positive Lebesgue measure. Moreover, the set of parameter values corresponding to stable periodic orbits is dense in any neighborhood of parameters with the chaotic attractor. For sufficiently large values of the Jacobian (very close to the classical H\'enon's parameters) an analogous result was obtained in \cite{GalTuc15} by means of computer-assisted proof methods. Thus, for specific parameter values, one can never be sure whether the chaotic attractor is observed or it just a transient chaos, after which orbits will tend to some stable periodic orbit with an extremely narrow absorbing domain.}, this attractor, together with the saddle point $O$, also contains the unstable manifold $W^u(O)$ and homoclinic points $h_i$ to $O$, i.e, such points where $W^u(O)$ intersects with the stable manifold $W^s(O)$. Figures~\ref{ris:Fig1}b and \ref{ris:Fig1}d show schematic representation of the homoclinic structure for orientable and nonorientable H\'enon attractors, respectively. Let us briefly describe their it.

The unstable invariant manifold $W^u(O)$ is divided by the point $O$ into two connected components -- separatrices $W^{u+}$ and $W^{u-}$. Since the unstable multiplier $\gamma$ of $O$ is negative ($\gamma<-1$), the separatrices $W^{u+}$ and $W^{u-}$ are invariant under $H^2$ and such that $H(W^{u+}) = W^{u-}$ and $H(W^{u-}) = W^{u+}$. This implies that points of $W^{u}$ jump under iterations of $H$ alternately from one separatrix to another.

The stable manifold $W^s(O)$ is also one-dimensional and it is divided by the point $O$ into two separatrices $W^{s+}$ and $W^{s-}$. For the orientable H\'enon attractor, the stable multiplier $\lambda$ is negative and thus, as for the unstable manifold, $H(W^{s+}) = W^{s-}$ and $H(W^{s-}) = W^{s+}$.

Let $h_1$ be an intersection point of $W^{u+}$ with $W^{s+}$. Then, $h_2$ is an intersection point of $W^{u-}$ with $W^{s-}$, since $H(W^{u+}) = W^{u-}$ and $H(W^{s+}) = W^{s-}$; $h_3$ is again an intersection point of $W^{u+}$ with $W^{s+}$, etc. Correspondingly, the points $h_1,h_2,...$ are homoclinic points of some homoclinic to $O$ orbit. Here, the points with odd indices $h_1, h_3, ...$ belong to the separatrix $W^{s+}$, while the points with even indices $h_2,h_4,...$ belong to $W^{s-}$, see Fig.~\ref{ris:Fig1}b. For the nonorientable H\'enon attractor, multiplier $\lambda$ is positive and, thus, all homoclinic points $h_1,h_2,...$ belong to one separatrix (e.g. $W^{s+}$) jumping from one unstable separatrix to another, see Fig.~\ref{ris:Fig1}d.

Similar homoclinic attractors are observed in many two-dimensional and multidimensional maps, as well as in Poincar\'e maps for various multidimensional systems of differential equations. Further, we will call them \textit{H\'enon-like attractors}. More precisely the H\'enon-like attractor can be defined as a homoclinic attractor containing the saddle fixed point (periodic orbit in the case of systems of differential equations) with a negative unstable multiplier and $\sigma < 1$.

\subsubsection{Phenomenological part of the scenario for hyperchaotic H\'enon-like attractor appearance}\label{sec:HenonPhen}

We start with the phenomenological part of the scenario leading to the appearance of hyperchaotic H\'enon-like attractor purposely skipping some details (accompanying bifurcations) which will be given in the framework of empirical part of the scenario in Sec.~\ref{sec:period_doubling}.

As for the H\'enon attractor, the first step in the framework of the scenario is the supercritical period-doubling bifurcation occurring with the stable fixed point $O$. Suppose that it happens at $\varepsilon = \varepsilon_{PD}$. After this bifurcation, the point $O$ becomes saddle of (2,1)-type and a stable period-2 orbit $(p_1, p_2)$ appears in its neighborhood, see Fig.~\ref{ris:Fig2}b. The saddle point $O$ has the following set of multipliers: $\gamma < - 1, -1 <\lambda_1 < 0$ and $|\lambda_2| < 1$. The unstable invariant manifold $W^u(O)$, separated by the point $O$ into two separatrices $W^{u+}$ and $W^{u-}$, is a segment with endpoints $p_1$ and $p_2$. Since this manifold corresponds to the negative multiplier $\gamma < -1$, we have a semi-local symmetry between the pair of separatrices $W^{u+}$ and $W^{u-}$ ($F(W^{u+}) = W^{u-}$ and $F(W^{u-}) = W^{u+}$). The stable invariant manifold $W^s$ corresponds to a pair of real multipliers $\lambda_1$ and $\lambda_2$. Suppose that $|\lambda_1| > |\lambda_2|$, i.e., $\lambda_1$ corresponds to the leading direction $W^{ls}$, and $\lambda_2$ -- to the strong stable manifold $W^{ss}$ tangent to the eigenvector corresponding to $\lambda_2$. Then, $\lambda_2 > 0$ if the map $F$ is orientable and $\lambda_2 < 0$, otherwise. Here, we consider only the case of orientation-preserving maps.

\begin{figure}[h!]
\begin{minipage}[h]{1\linewidth}
\center{\includegraphics[width=1\linewidth]{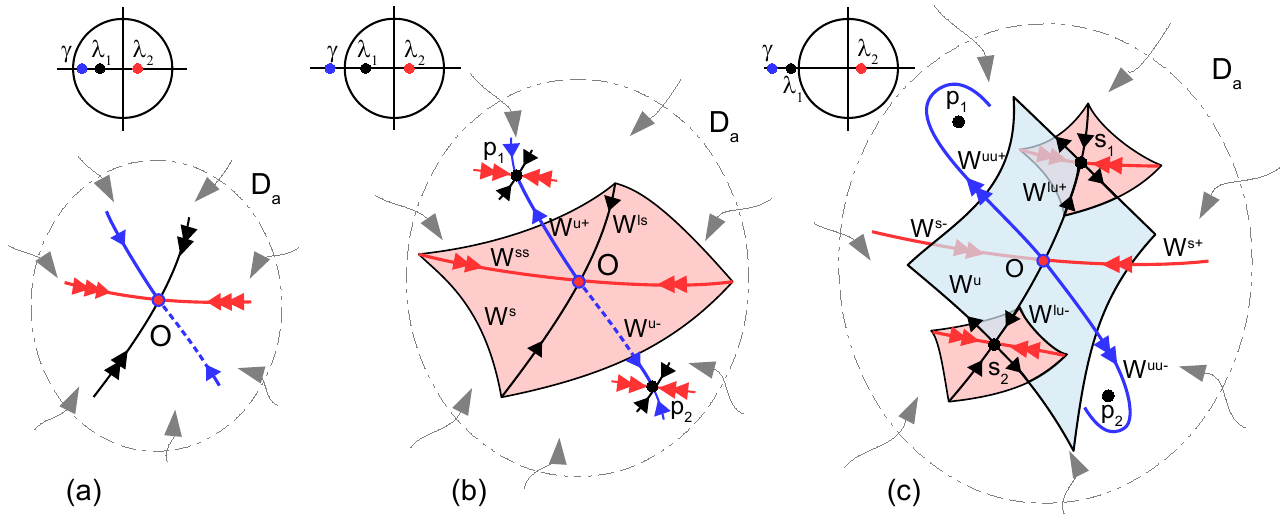}}
\end{minipage}
\caption{\footnotesize Illustrations to the phenomenological part of the scenario of hyperchaotic H\'enon-like attractor appearance. (a) $O$ is a stable fixed point, it resides in some absorbing domain $D_a$; (b) period-doubling bifurcation with the point $O$ occurs: a stable period-2 orbit $(p_1,p_2)$ is born, while the point $O$ becomes saddle with a negative unstable multiplier; (c) period-doubling bifurcation with the saddle point $O$ occurs: a saddle period-2 orbit $(s_1,s_2)$ is born, while the point $O$ becomes saddle with a pair of negative unstable multipliers.}
\label{ris:Fig2}
\end{figure}

\begin{remark}
The first step in the framework of this scenario is the same as for the creation of the so-called discrete Lorenz and figure-eight homoclinic attractors \cite{GOST05, GGKT14, GG16}. Both these attractors are remarkable since they can be pseudohyperbolic \cite{TS98, TS08}. Pseudohyperbolicity is a weak version of hyperbolicity. Chaotic dynamics of pseudohyperbolic attractors persist under small perturbations (as for hyperbolic attractors), despite the possible occurrence of homoclinic tangencies inside them. Note that the saddle fixed point belonging to both these attractors should be area-expanding, which means that $\gamma\lambda_1 > 1$.
\end{remark}

The next principal bifurcation in the framework of the scenario is a second supercritical period-doubling bifurcation with the saddle fixed point $O$. Suppose that it occurs at $\varepsilon = \varepsilon^{PD}$. After this bifurcation, the point $O$ becomes saddle of (1,2)-type a period-2 saddle orbit $(s_1, s_2)$ of (2,1)-type appears in its neighborhood, see Fig.~\ref{ris:Fig2}?. For $\varepsilon > \varepsilon^{PD}$, both unstable multipliers $\gamma$ and $\lambda_1$ of the point $O$ are negative, while the stable multiplier $\lambda_2$ is positive since we consider the orientable case. The following conditions on the multipliers are met here
$$
\gamma < \lambda_1 < -1, \; 0 < \lambda_2 <1.
$$

A restriction $F_u$ of the initial map $F$ into the local unstable manifold $W^u_{loc}(O)$ has a fixed point $\tilde O = O\cap W^u_{loc}$ which is the unstable node with a pair of multipliers $\gamma < \lambda_1 < -1$. Thus, in $W^u_{loc}$, there are a strong unstable invariant manifold $W^{uu}$, which is tangent to the eigenvector corresponding to the unstable multiplier $\gamma$, and a leading unstable direction $W^{lu}$ corresponding to the multiplier $\lambda_1$. Also note that the curves $W^{uu}$ and $W^{lu}$ divide $W^{u}_{loc}$ into four fragments $\Pi_1$, $\Pi_2$, $\Pi_3$, and $\Pi_4$ and, since both multipliers of $\tilde O$ are negative, $F_u(\Pi_1) = \Pi_3$, $F_u(\Pi_3) = \Pi_1$, $F_u(\Pi_2) = \Pi_4$, and $F_u(\Pi_4) = \Pi_2$. All orbits in $W^u_{loc}$, except those that belong to $W^{uu}$, tend to the node $\tilde O$ (in backward time) along the smooth cubic parabola-like curves. These curves are tangent to the leading direction $W^{lu}$ of $\tilde O$, see Fig.~\ref{ris:Fig3}a. If $h_1 \in \Pi_1$ is one of such points, then each its odd (even) iteration under $F^{-1}_u$ tends to $\tilde O$ along right (left) branch of this parabola-like curve staying in $\Pi_3$ ($\Pi_1$). It is important to note that since, by the moment, only orbits belonging to $W^{s}(O)$ tend to the saddle fixed point $O$, for other orbits there is no mechanism to return into the neighborhood of this point.

\begin{figure}[h!]
\begin{minipage}[h]{1\linewidth}
\center{\includegraphics[width=0.7\linewidth]{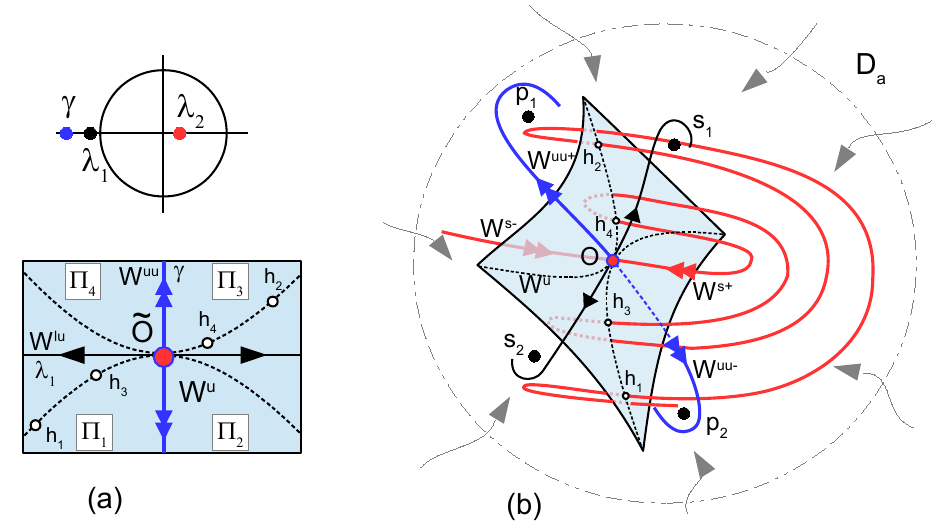}}
\end{minipage}
\caption{\footnotesize Continuation of Fig.~\ref{ris:Fig2}: (a) schematic representation of the local unstable manifold $W^u_{loc}(O)$ on which the point $\tilde O$ is an unstable node; (b) schematic representation of the homoclinic structure for the hyperchaotic H\'enon-like attractor.}
\label{ris:Fig3}
\end{figure}

We suppose that such mechanism appears at $\varepsilon = \varepsilon^{H}$ due the emergence of a homoclinic orbit to $O$, see Fig.~\ref{ris:Fig3}b, when the stable manifold $W^{s}(O)$, separated by the point $O$ into two separatrices $W^{s+}$ and $W^{s-}$, starts to intersect with the unstable manifold $W^{u}(O)$. Suppose that $h_1$ is the first intersection point of $W^{s+}$ with $W^{u}$ at $\Pi_1$, see Fig.~\ref{ris:Fig3}a. Then (since $F^{-1}(W^{s+}) = W^{s+}$ and $F_u(\Pi_1) = \Pi_3$), $h_2 = F^{-1}(h_1)$ is the intersection point of the same stable separatrix $W^{s+}$ with $W^{u}$ in $\Pi_3$; $h_3 = F^{-1}(h_2)$ is again an intersection point of $W^{s+}$ with $W^{u}$ in $\Pi_1$, etc. Correspondingly, the points $h_1,h_2,...$ are homoclinic points of some homoclinic to $O$ orbit, see Fig.~\ref{ris:Fig3}b.

At $\varepsilon > \varepsilon^{H}$, this homoclinic orbit gives a nontrivial hyperbolic subset of (1,2)-type which, we suppose, becomes attractive, i.e., we observe, in the absorbing domain $D_a(O)$, the homoclinic attractor containing the saddle point $O$ with a pair of negative unstable multipliers. We call attractors of such type \textit{hyperchaotic H\'enon-like attractors}. If the numerically obtained orbits in the attractor provide a sufficiently long time near the hyperbolic subset of (1,2)-type, then a pair of Lyapunov exponents becomes positive, i.e., we observe hyperchaos in numerical experiments\footnote{We note that the described hyperchaotic attractor typically is not hyperbolic. Together with periodic saddle orbits of (1,2)-type it also contains periodic saddle orbits of a (2,1)-type and homoclinic (heteroclinic) tangencies between invariant manifolds of various saddles. Moreover, it is possible to show that such attractors also contain heterodimensional cycles connecting these saddles \cite{Li16, LiT17, LiT20}. Such cycles robustly persist under small perturbation of the system guaranteeing the presence of a nontrivial hyperbolic subset of (1,2)-type.}.

In the following section, we show how hyperchaotic H\'enon-like attractors can naturally appear in multidimensional systems demonstrating cascades of period-doubling bifurcations. In Sec.~\ref{sec:Shil_emp}, we generalize this scenario to attractors developing from stable periodic orbits which occur inside Arnold tongues. In Section~\ref{sec:SmallB}, we demonstrate its implementation for periodic resonant orbits occurring inside Arnold tongues of the three-dimensional Mir\'a map \eqref{eq:MiraMap}.

\subsubsection{Empirical part of the scenario for hyperchaotic H\'enon-like attractor appearance} \label{sec:period_doubling}

\begin{figure}[h!]
\begin{minipage}[h]{1\linewidth}
\center{\includegraphics[width=1\linewidth]{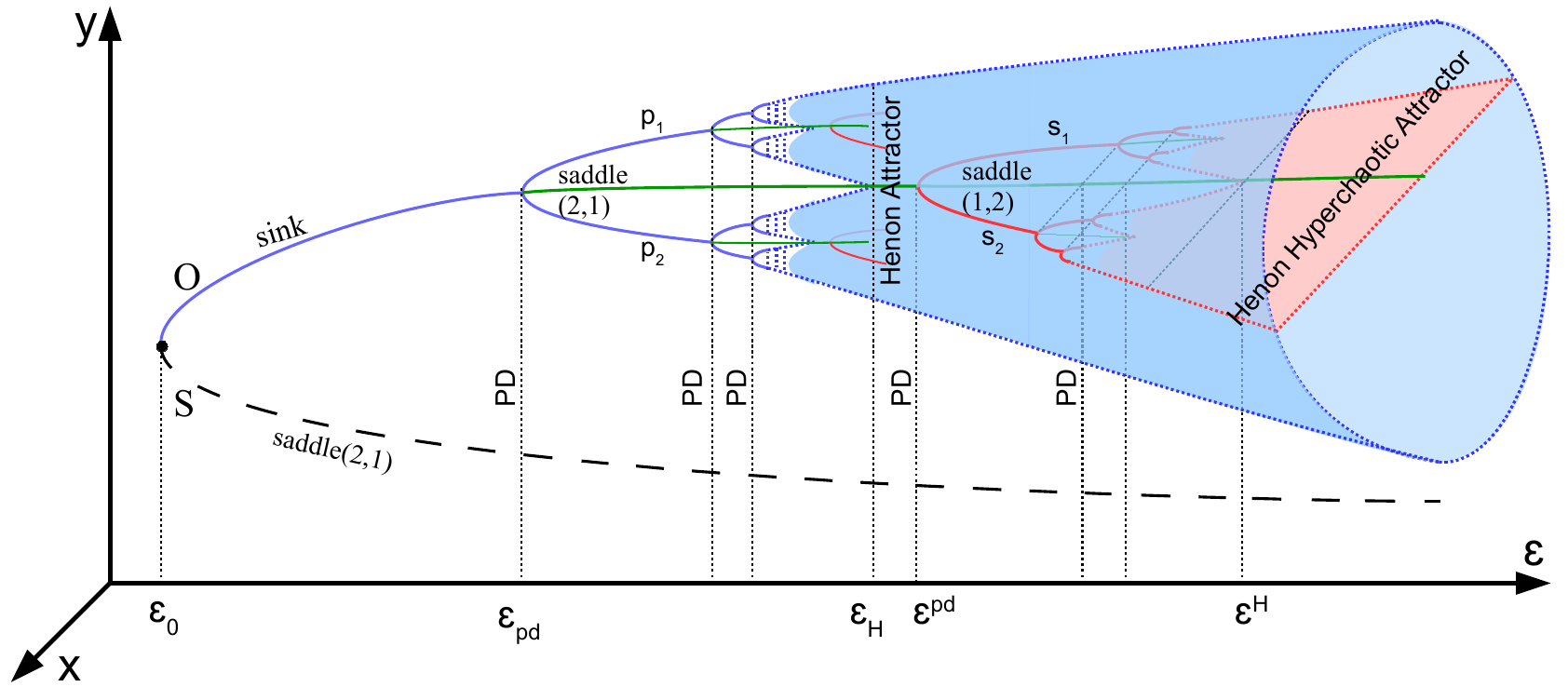}}
\end{minipage}
\caption{\footnotesize Schematic representation of the scenario of hyperchaotic H\'enon-like attractor appearance as a result of a cascade of period-doubling bifurcations with the stable fixed point $O$ followed by cascades of period-doubling bifurcations with periodic saddle orbits of (2,1)-type.}
\label{ris:Fig4}
\end{figure}

Here we show how the attractor described above can appear in multidimensional systems demonstrating transition to chaos via cascades of period-doubling bifurcations. Suppose that, at $\varepsilon = \varepsilon_H$, a H\'enon-like attractor appears after a cascade of period-doubling bifurcations followed by a cascade of heteroclinic band-merging bifurcations, see Fig.~\ref{ris:Fig4}. Note that for one-dimensional maps the exponential convergence of these bifurcations is well-known fact \cite{lanford1982computer, campanino1981existence, eckmann1987complete}. For two-dimensional maps such universality is proved in Refs.~\cite{CarLyubMart05, LyubMart11}. As we know, for three-dimensional maps there are no rigorous results on this topic. However it is believed that the exponential convergence with the same universal property is also observed in this case \cite{collet1981period}.

\begin{figure}[h!]
\begin{minipage}[h]{1\linewidth}
\center{\includegraphics[width=1\linewidth]{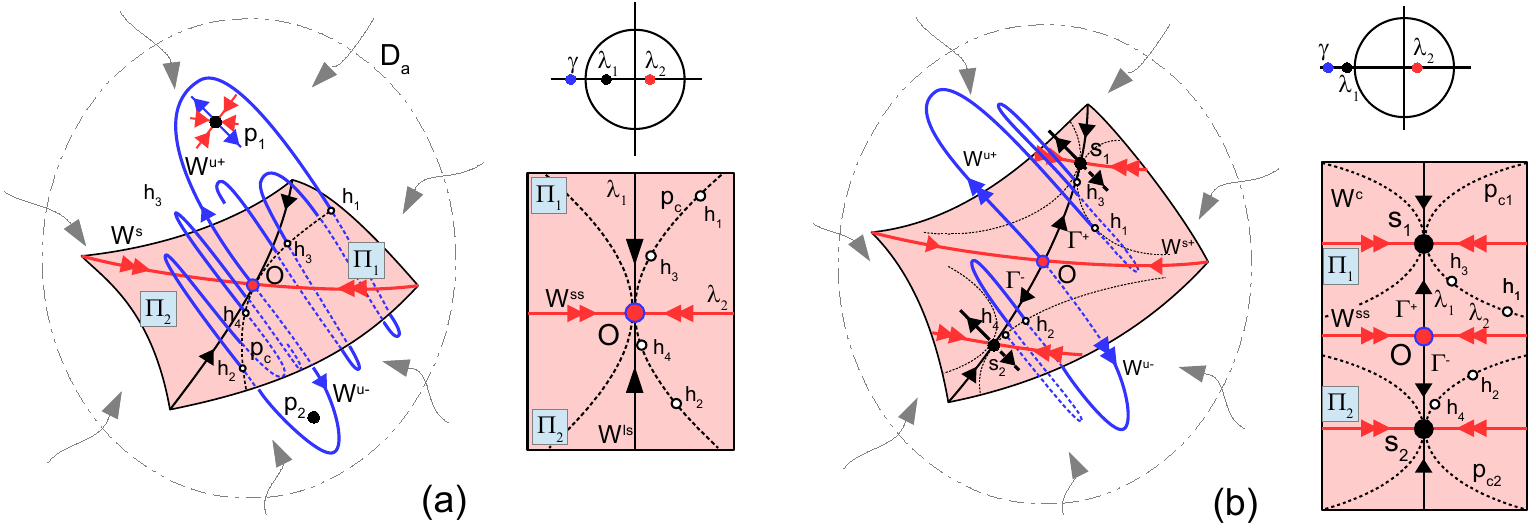}}
\end{minipage}
\caption{\footnotesize (a) Schematic representation of the homoclinic structure for a H\'enon-like attractor of a three-dimensional map; (b) invariant manifolds after the period-doubling bifurcation with the saddle fixed point $O$.}
\label{ris:Fig5}
\end{figure}

A homoclinic structure for the H\'enon-like attractor is shown schematically in Fig.~\ref{ris:Fig5}a (left panel). Here, as in the two-dimensional case, the unstable multiplier $\gamma$ of the fixed point $O$ is negative ($\gamma < -1$) and the stable ones $\lambda_1$ and $\lambda_2$ are real. Depending on the signs of stable multipliers, homoclinic structures for the H\'enon-like attractor can be of four possible types (two in orientable and two in nonorientable cases). Here we consider only one orientable case characterizing by the following values of multipliers $\gamma < -1 < \lambda_1 < 0 < \lambda_2 < 1, |\lambda_2| < |\lambda_1|$ (see the schematic location of the multipliers at the top-right insert of Fig.~\ref{ris:Fig5}a). This attractor is a three-dimensional generalization of the classical orientable H\'enon attractor.

Here the stable invariant manifold $W^s(O)$ is two-dimensional. A restricted map $F_s$ of the initial map $F$ into $W^s_{loc}(O)$ has a fixed point $\tilde O = O\cap W^s_{loc}$ which is a stable node with a pair of multipliers $\lambda_1$ and $\lambda_2$. According to the condition $|\lambda_1| > |\lambda_2|$, in $W^s_{loc}$ there are the strong stable invariant manifold $W^{ss}$, which is tangent to the eigenvector corresponding to $\lambda_2>0$, and the leading direction $W^{ls}$ corresponding to $\lambda_1 < 0$. The curve $W^{ss}$ divides $W^s_{loc}$ into two parts $\Pi_1$ and $\Pi_2$, and, since the leading multiplier is negative ($\lambda_1 < 0$), $F_s(\Pi_1) = \Pi_2$ and $F_s(\Pi_2) = \Pi_1$. All orbits in $W^s_{loc}$, except those that belong to $W^{ss}$ tend to the node $\tilde O$ along the parabola-like curves which are tangent to the leading direction $W^{ls}$, see the right panel in Fig.~\ref{ris:Fig5}a. If we take some point $h_1$ belonging to one such curves $p_c$ and consider its images under $F_s$, we obtain an orbit $h_1, h_2, ...$ tending to $\tilde O$ on one (e.g. right) side of $W^{ls}$ and jumping from one branch of $p_c$ to another after each iteration. Suppose, $h_1$ is the intersection point of $W^{u+}$ with $W^{s}$ at its upper part $\Pi_1$. Then, (since $F(W^{u+}) = W^{u-}$ and $F_s(\Pi_1) = \Pi_2$) $h_2 = F(h_1)$ is the intersection point of $W^{u-}$ with $W^{s}$ in its bottom part $\Pi_2$; $h_3 = F(h_2)$ is again an intersection point of $W^{u+}$ with $W^{s+}$ in $\Pi_1$, etc.

Then, we suppose that the saddle point $O$ of (2,1)-type undergoes a cascade of period-doubling bifurcations, see the schematic phase portrait after the first period-doubling bifurcation in Fig.~\ref{ris:Fig5}b. After each such bifurcation, the periodic orbit gets the two-dimensional unstable manifold and a double-period orbit of (2,1)-type appears in its neighborhood. Then, we suppose that the period-$2^n$ saddle orbits of (2,1)-type emerging after the corresponding period-doubling bifurcations, as well as other saddles of (2,1)-type inside the attractor, also undergo full cascades of period-doubling bifurcation after which all these orbits become of (1,2)-type. It is natural to assume that these cascades have the same universal properties of exponential convergence as the cascades for stable periodic orbits \cite{collet1981period}. We assume that after it, we obtain a hyperchaotic attractor which consists of (infinitely) many components separated by isolated periodic saddles of (1,2)-type (similar to chaotic attractor in the H\'enon map before the cascade of heteroclinic band-merging bifurcations).

Further, we suppose that this multi-component attractor goes through a cascade of heteroclinic band-merging bifurcations. These bifurcations are similar, in a sense, to those ones which lead to the pairwise merger of components of a chaotic attractor arising after the successive cascade of period-doubling bifurcations within the second part of the scenario resulting in the birth of a H\'enon-like attractor. After each such bifurcation, the number of components decreases by a factor of two.

The final bifurcation here, resulting in the merger of two last components, leads to the emergence of homoclinic intersections between $W^s(O)$ and $W^u(O)$. As a result, the hyperchaotic H\'enon-like attractor occurs at $\varepsilon = \varepsilon^{H}$. Note, that this intersection persists at some parameter region. However, the hyperchaotic H\'enon-like attractor (as the classical H\'enon attractor) is not pseudohyperbolic (robustly chaotic). It contains the saddle fixed point $O$ with the two-dimensional unstable invariant manifold and (as we have an attractor) with such a saddle index $\rho = |\lambda_1 \lambda_2|$, that $\rho < 1$. As was shown in \cite{GST96, GST08, GST97}, stable periodic orbits can appear in this case under arbitrarily small perturbations due to bifurcations of homoclinic tangencies.

\subsection{Discrete Shilnikov attractor}\label{sec:ShilAttrScen}

As shown at the beginning of this section, a period-doubling bifurcation can be the first step in the framework of onset of H\'enon-like and hyperchaotic H\'enon-like attractors. A natural question arises here. Which homoclinic attractors can appear if the first codimension-one bifurcation is a Neimark-Sacker bifurcation? The answer to this question was given in \cite{GGS12} (see also \cite{GGKT14}, and \cite{GG16}) where a scenario for the appearance of the so-called \textit{discrete Shilnikov attractor} containing a saddle-focus fixed point of (1,2)-type was proposed. This scenario goes back to the paper by Shilnikov \cite{Sh86} where a similar scenario was proposed for one-parameter families of three-dimensional flows. Since this scenario plays an important role, in the development of hyperchaos \cite{GSKK19}, let us briefly describe it.

\subsubsection{Phenomenological part of the scenario of discrete Shilnikov attractor appearance} \label{sec:Shil_phen}

\begin{figure}[h!]
\begin{minipage}[h]{1\linewidth}
\center{\includegraphics[width=1\linewidth]{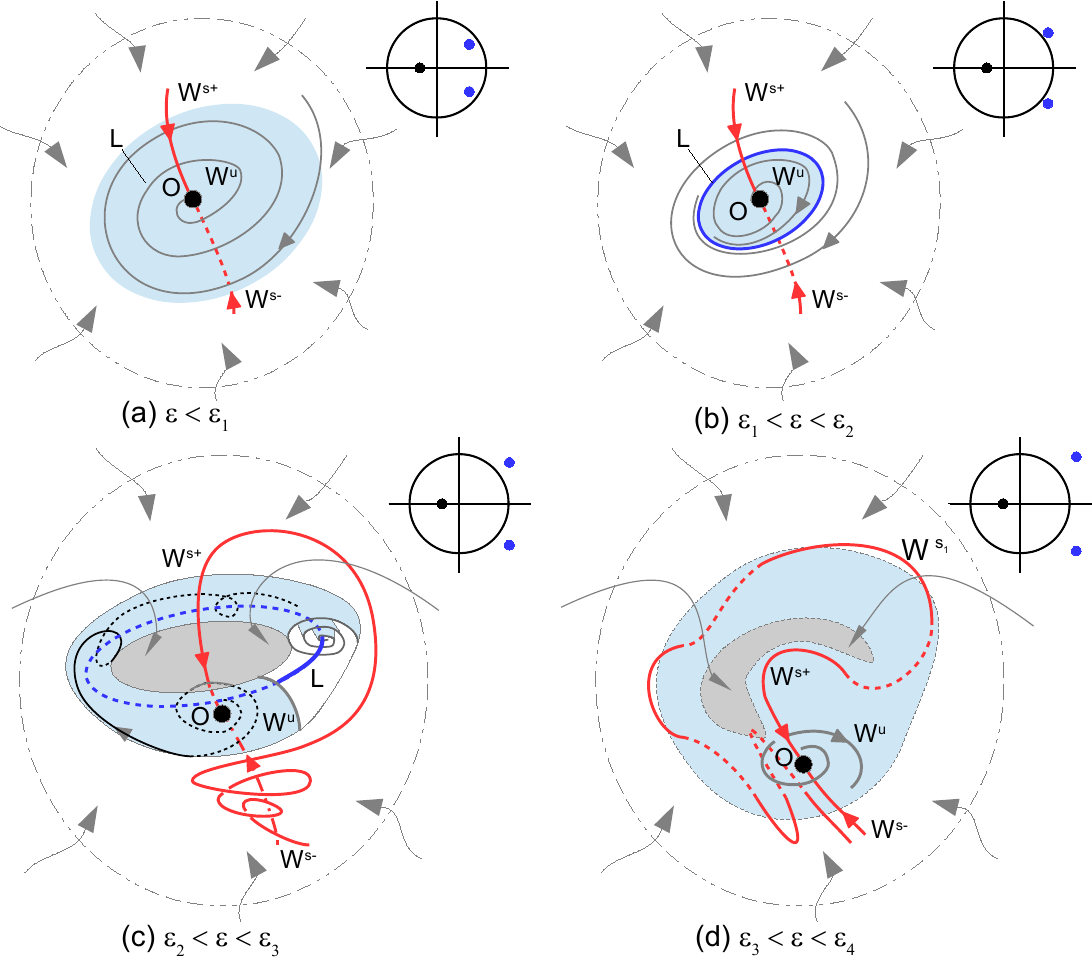}}
\end{minipage}
\caption{\footnotesize Illustration towards the scenario of the appearance of a discrete Shilnikov attractor containing a saddle-focus fixed point of (1,2)-type. This figure is taken from \cite{GG16}.}
\label{ris:Fig6}
\end{figure}

Suppose that at $\varepsilon < \varepsilon_1$ the fixed point $O$ is stable but focal (a pair of its multipliers is complex-conjugate), see Fig.~\ref{ris:Fig6}a. At $\varepsilon = \varepsilon_1$ it undergoes the supercritical Neimark-Sacker bifurcation. As a result, this fixed point becomes a saddle-focus of (1,2)-type, and a stable invariant curve $L$ is born in its neighborhood, see Fig.~\ref{ris:Fig6}b. Note that after the birth, this curve has a nodal type: the two-dimensional unstable invariant manifold $W^u(O)$ is a disc with an edge on $L$. Then, we suppose that at $\varepsilon = \varepsilon_2$ the stable curve becomes focal, and, as a result, $W^u(O)$ starts to wind on it forming the so-called ``Shilnikov whirlpool'', see Fig.~\ref{ris:Fig6}c. All orbits from the absorbing domain $D_a(O)$ (except the stable separatrix $W^{s-}(O)$) are tightened by this whirlpool. With further increase in $\varepsilon$, the size of whirlpool is increased, and finally, at $\varepsilon = \varepsilon_3$, the stable separatrix $W^{s+}(O)$ touches $W^u(O)$. As a result, at some interval $\varepsilon_3 < \varepsilon < \varepsilon_4$, the fixed point $O$ has a transversal homoclinic structure. By Smale and Shilnikov \cite{Sm67, Sh67}, in the neighborhood of the transversal intersection $W^u(O) \cap W^{s+}(O)$ there exist countably many periodic orbits of the same type as the fixed point $O$ which is a saddle-focus of (1,2)-type.

If the stable invariant curve $L$ breaks down (by Afraimovich-Shilnikov \cite{AfrSh83b} or due to some other scenario \cite{Komuro16, ZM08, ZM09, ZM10}) giving chaotic attractor, then, on some subinterval inside $\varepsilon \in (\varepsilon_3, \varepsilon_4)$, this attractor can contain the fixed point $O$ together with the nontrivial hyperbolic set of (1,2)-type. Therefore, this discrete Shilnikov attractor can be hyperchaotic.

\begin{remark}
In the framework of the described scenario we suppose that the fixed point $O$ undergoes the supercritical Neimark-Sacker bifurcation. However, it can undergo the subcritical Neimark-Sacker bifurcation, i.e., the point $O$ can sharply lose the stability due to the merger with a saddle invariant curve existing in a neighborhood of $O$. In this case, a discrete Shilnikov attractor can appear suddenly. A similar scenario (but for a three-dimensional flow system) was observed e.g. in \cite{KKL19}.
\end{remark}

\subsubsection{Empirical part of the scenarios for hyperchaotic and flow-like discrete Shilnikov attractors appearance} \label{sec:Shil_emp}

\begin{figure}[h!]
\begin{minipage}[h]{1\linewidth}
\center{\includegraphics[width=1\linewidth]{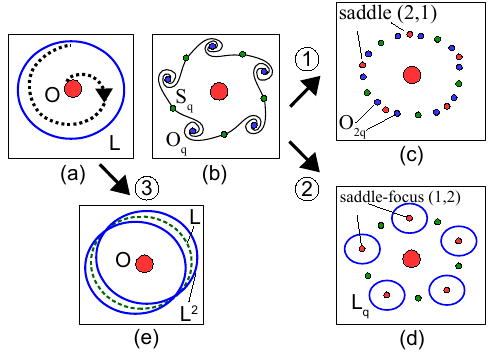}}
\end{minipage}
\caption{\footnotesize Schematic representation of two possible mechanisms of the destruction of the invariant curve $L$: (a) $L$ is a nonresonant (ergodic) invariant curve; (b) $L$ is a resonant curve (a pair of stable and saddle period-$q$ orbits $O_q$ and $S_q$ appears on it via the saddle-node bifurcation); (c) $O_q$ undergoes the supercritical period-doubling bifurcation (d) $O_q$ undergoes the supercritical Neimark-Sacker bifurcation; (e) $L$ undergoes a supercritical period-doubling (length-doubling) bifurcation.}
\label{ris:Fig7}
\end{figure}

The empirical part of the scenario for discrete Shilnikov attractors appearance describes specific sequences of bifurcations responsible for chaotization of the invariant curve $L$ (Fig.~\ref{ris:Fig7}a) and the absorption of the saddle-focus fixed point $O$ by the resulting attractor. Let us describe possible scenarios leading to the emergence of (i) hyperchaotic and (ii) flow-like strange attractors on the base of $L$. Suppose that before the destruction, an invariant curve becomes resonant: a pair of stable and saddle period-$q$ orbits $O_q$ and $S_q$ appears on it as a result of the saddle-node bifurcation. In this case, the curve is formed by the the closure of the unstable one-dimensional invariant manifold of $S_q$, see Fig.~\ref{ris:Fig7}b. Note that the stable two-dimensional invariant manifold of $S_q$ forms the boundary of the absorbing domain for the stable point $O_q$ or attractors developing from it in accordance with one of two following scenarios.
\begin{enumerate}
\item The point $O_q$ undergoes the supercritical period-doubling bifurcation after which it becomes saddle of (2,1)-type and a stable period-$2q$ orbit $O_{2q}$ is born in its neighborhood, see Fig.~\ref{ris:Fig7}c, i.e., in this case we have a hypothetical possibility to obtain the hyperchaotic H\'enon-like attractor on the base of $O_q$ (see Sec.~\ref{sec:period_doubling} for more detail). However, we note that the cascade of period-doubling bifurcations with a stable periodic orbit developed from $O_q$ can terminate at some moment by a Neimark-Sacker bifurcation after which the resulting periodic orbit gets the two-dimensional unstable manifold immediately, and a stable multi-component invariant curve is born. In its turn, this curve, before the destruction, can also become resonant, for it one of two cases under consideration is applicable, and so on.
\item The point $O_q$ undergoes the supercritical Neimark-Sacker bifurcation after which it becomes saddle-focus of (1,2)-type and a stable $q$-component invariant curve $L_q$ is born in its neighborhood, see Fig.~\ref{ris:Fig7}d. So, here we have possibility to obtain a $q$-component Shilnikov attractor. It is important to note, that the curve $L_q$, before the destruction can also become resonant, for it one of two cases under consideration is applicable, and so on.
\end{enumerate}
Concerning the resonant saddle orbit $S_q$ (as well as other resonant saddle orbits arising in the framework of the described above sequence), it undergoes the cascade of period-doubling bifurcations resulting in a one-component nontrivial hyperbolic set SH$_q$(1,2) of (1,2)-type (the same as with the saddle fixed point in the framework of the scenario described in Sec.~\ref{sec:period_doubling}). Touching the one-dimensional stable manifolds of SH$_q$(1,2) the attractor developed from $O_q$ undergoes crisis and collides into a one-component attractor which, after that, contains SH$_q$(1,2).

\begin{enumerate}
\item[3. ] Note that before becoming resonant the curve $L$ can go through a long sequence of period-doubling (length-doubling) bifurcations. After each such bifurcation the invariant curve becomes saddle, and a stable double-round invariant curve appears in its neighborhood (Fig.~\ref{ris:Fig7}e). Such transition to chaotic attractors (including Shilnikov ones), is observed quite often (see e.g. \cite{ACS85, Kaneko83, AnNik05, BSV05, BSV10, BKS16, St2, grines2022origin}) and can lead to the birth of flow-like chaotic attractors possessing one positive and one zero (indistinguishable from zero in numerical experiments) Lyapunov exponents (for ODE systems strange attractors with two zero Lyapunov exponent is born via this scenario). In Sec.~\ref{sec:VarShil} we demonstrate and study this phenomenon in more detail.
\end{enumerate}

\section{Three-dimensional Mir\'a map: main bifurcations and dynamical regimes}\label{sec:MiraMapBiff}

In this section we study main bifurcations and discuss the most interesting dynamical regimes in the three-dimensional Mir\'a map \eqref{eq:MiraMap}. Recall, that this map has a constant Jacobian, $J=B$. Various types of attractors are possible only when the map is dissipative, i.e., when $|B| < 1$. In this paper we consider only orientation preserving maps, when $B$ is positive ($0 < B < 1$).

Map \eqref{eq:MiraMap} has up to two fixed points $O_+(x_+,y_+,z_+)$ and $O_-(x_-,y_-,z_-)$ with coordinates
\begin{equation}
    x_{\pm}=y_{\pm}=z_{\pm}=\frac{M_2 + B - 1}{2} \pm \sqrt{\frac{(1 - B - M_2)^2}{4} + M_1}.
\end{equation}
As one can see, these points exist only when the root expression is non-negative. Points $O_{\pm}$ are born under the saddle-node (tangent) bifurcation (when one multiplier is equal to $1$) on the surface
\begin{equation}
    \mbox{SN}: M_1 = - \frac{(1 - B - M_2)^2}{4}.
    \label{eq:SN}
\end{equation}
A period-doubling bifurcation occurs when a multiplier is equal to $-1$ on the surface
\begin{equation}
    \mbox{PD}: M_1 = \frac{3(B + M_2)^2 + 2(B + M_2) - 1}{4}.
    \label{eq:PD}
\end{equation}
The third codimension-one bifurcation appears when a pair of multipliers becomes equal to $e^{\pm i \phi}, \phi \in (0, \pi)$. This case corresponds to a Neimark-Sacker bifurcation which appears on the surface
\begin{equation}
    \mbox{NS}: M_1 = \frac{(2 - M_2 - B + M_2B - B^2)^2 - (1 - B - M_2)^2}{4}, \quad |M_2 - B| < 2.
    \label{eq:NS}
\end{equation}

\begin{figure}[h!]
\begin{minipage}[h]{1\linewidth}
\center{\includegraphics[width=1\linewidth]{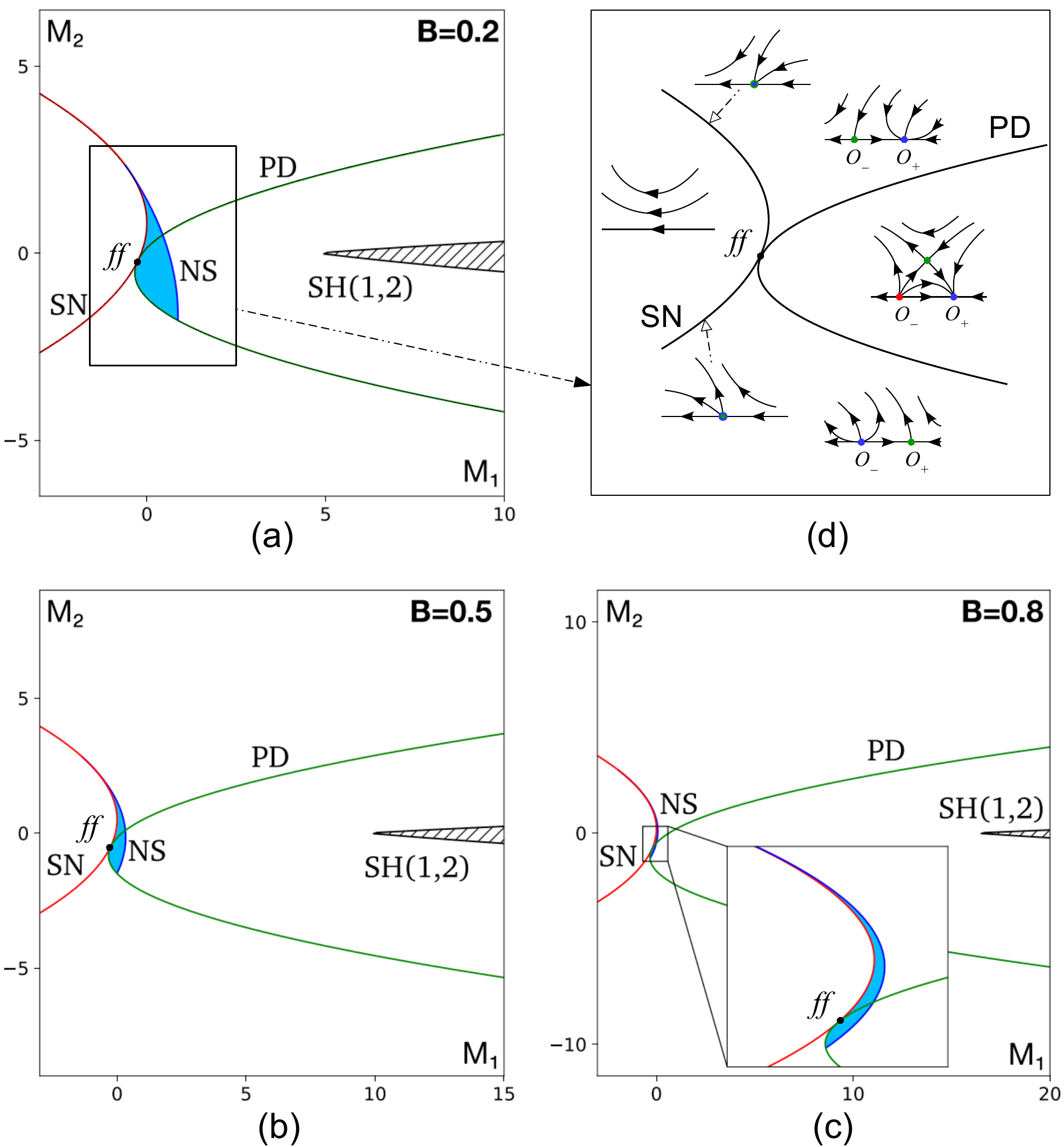}}
\end{minipage}
\caption{\footnotesize Some bifurcation curves for map \eqref{eq:MiraMap} for various values of parameter $B$: (a) $B=0.2$, (b) $B=0.5$, (c) $B=0.8$. SN, PD, and NS are the saddle-node, period-doubling, and Neimark-Sacker bifurcation curves. The region of stability of the point $O_+$ is colored in blue. The region SH(1,2) corresponds to the existence of a nontrivial hyperbolic subset of (1,2)-type. (d) Sketch of the bifurcation diagram near the fold-flip bifurcation point $ff$. Here, the behavior of orbits on the corresponding center 2D manifold is shown; this invariant manifold is asymptotically stable in our case.}
\label{ris:Fig8}
\end{figure}

Note that relations \eqref{eq:SN}-\eqref{eq:NS} define in the three-dimensional parameter space $(M_1, M_2, B)$ the region of stability of the point $O_+$. Figure~\ref{ris:Fig8} shows several slices of this space for various values of the parameter $B$. In this figure the region of stability for $O_+$ is colored in blue, let us explain its boundaries.

The curves SN and PD are tangent at the codimension-two point $ff$ corresponding to a \textit{fold-flip bifurcation}. As is known from the paper \cite{kuznetsov2004fold} (see also \cite{gheiner1994codimension}), depending on coefficients of the corresponding normal form, this bifurcation can be of four possible types. Our studies show that here we observe the fourth case from \cite{kuznetsov2004fold}. The sketch for this bifurcation is shown in Fig.~\ref{ris:Fig8}b. When crossing the SN-curve above the point $ff$, the stable fixed point $O_+$ is born together with the saddle $O_-$. The point $O_-$ undergoes a period-doubling bifurcation on the curve PD above $ff$, as a result this point becomes of saddle of (1,2)-type and a period-2 saddle orbit of type (2,1) is born in its neighborhood. This period-2 orbit merges with the stable fixed point $O_+$ on the bottom piece of the PD-curve (below $ff$), after which the point $O_+$ becomes saddle of (2,1)-type. Finally, the saddle points $O_+$ and $O_-$ merges on the bottom piece of the SN-curve and both disappear to the left of this curve.

In Fig.~\ref{ris:Fig8} we also plot the region SH(1,2) inside which the nonwandering set of map \eqref{eq:MiraMap} is a hyperchaotic hyperbolic set consisting of saddle periodic orbits of (1,2)-type. We also call this set Smale horseshoe of (1,2)-type. According to the theorem by Gonchenko and Li (th.~2 in \cite{GonLi10}), this region is bounded by the surface
$$
M_1 = \left(\tilde{\rho} + \sqrt{\tilde{\rho}^2 + \frac{1}{4}}\right)^2 - (1 - M_2 - B)\left(\tilde{\rho} + \sqrt{\tilde{\rho}^2 + \frac{1}{4}}\right),
$$
where
$$
\tilde{\rho} = \frac{3 + 5(|B| + |M_2|)}{3 + 4(|B| + |M_2|)}(1 + |B| + |M_2|).
$$
In this paper we study evolution of attractors along pathways from the region of stability of the fixed point $O_+$ towards the region SH(1,2).

As shown below, all interesting dynamical regimes in map~\eqref{eq:MiraMap} for $0<B<1$ are associated with bifurcations of the fixed point $O_+$. Thus, it is convenient to shift this point to the origin which gives the following representation for map \eqref{eq:MiraMap}

\begin{equation}
    \begin{cases}
        \bar{x} = y \\
        \bar{y} = z \\
        \bar{z} = Bx + Cy + Az - y^2.
    \end{cases}
    \label{eq:Henon3DMap}
\end{equation}

\begin{figure}[h!]
\begin{minipage}[h]{1\linewidth}
\center{\includegraphics[width=0.6\linewidth]{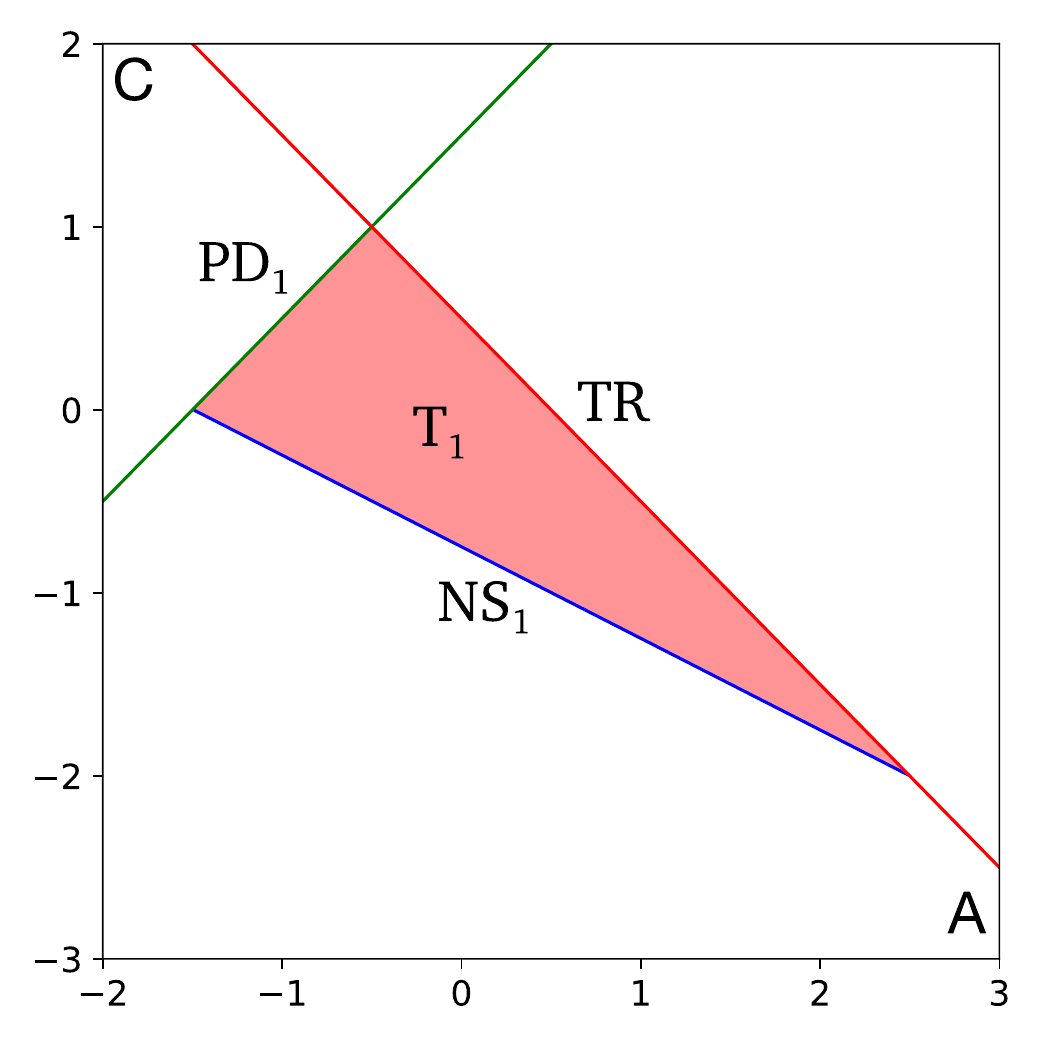}}
\end{minipage}
\caption{\footnotesize Some bifurcation curves of map \eqref{eq:Henon3DMap}, $B=0.5$. The curves TR, PD$_1$, and NS$_1$ correspond to the transcritical saddle-node, period-doubling, and supercritical Neimark-Sacker bifurcations. The region of stability of the point $O_+$ is colored in pink.}
\label{ris:Fig9}
\end{figure}

The main difference between maps \eqref{eq:MiraMap} and \eqref{eq:Henon3DMap} is that in the last map both fixed points $O_+$ and $O_-$ always exist. The point $O_+$ becomes stable here under a transcritical saddle-node bifurcation (but note via a saddle-node bifurcations as for map \eqref{eq:MiraMap}) on a plane
\begin{equation}
    \mbox{TR}: C = 1 - A - B.
    \label{eq:TR}
\end{equation}
As in map \eqref{eq:MiraMap}, other boundaries of stability of this point are determined by a period-doubling bifurcation on a plane
\begin{equation}
    \mbox{PD}_1: C = 1 + A + B
    \label{eq:PD1}
\end{equation}
and a Neimark-Sacker bifurcation on a surface
\begin{equation}
    \mbox{NS}_1: C = B^2 - AB - 1, \quad -2<A-B<2,
    \label{eq:NS1}
\end{equation}
see bifurcation diagram for $B=0.5$ in Fig.~\ref{ris:Fig9}. Respectively, relations \eqref{eq:TR}--\eqref{eq:NS1} define in the three-dimensional parameter space $(A, B, C)$ the region of stability of the fixed point $O_+$ bounded by the transcritical saddle-node (TR), period-doubling (PD$_1$), and Neimark-Sacker (NS$_1$) bifurcation curves.

\begin{figure}[h!]
\begin{minipage}[h]{1\linewidth}
\center{\includegraphics[width=1\linewidth]{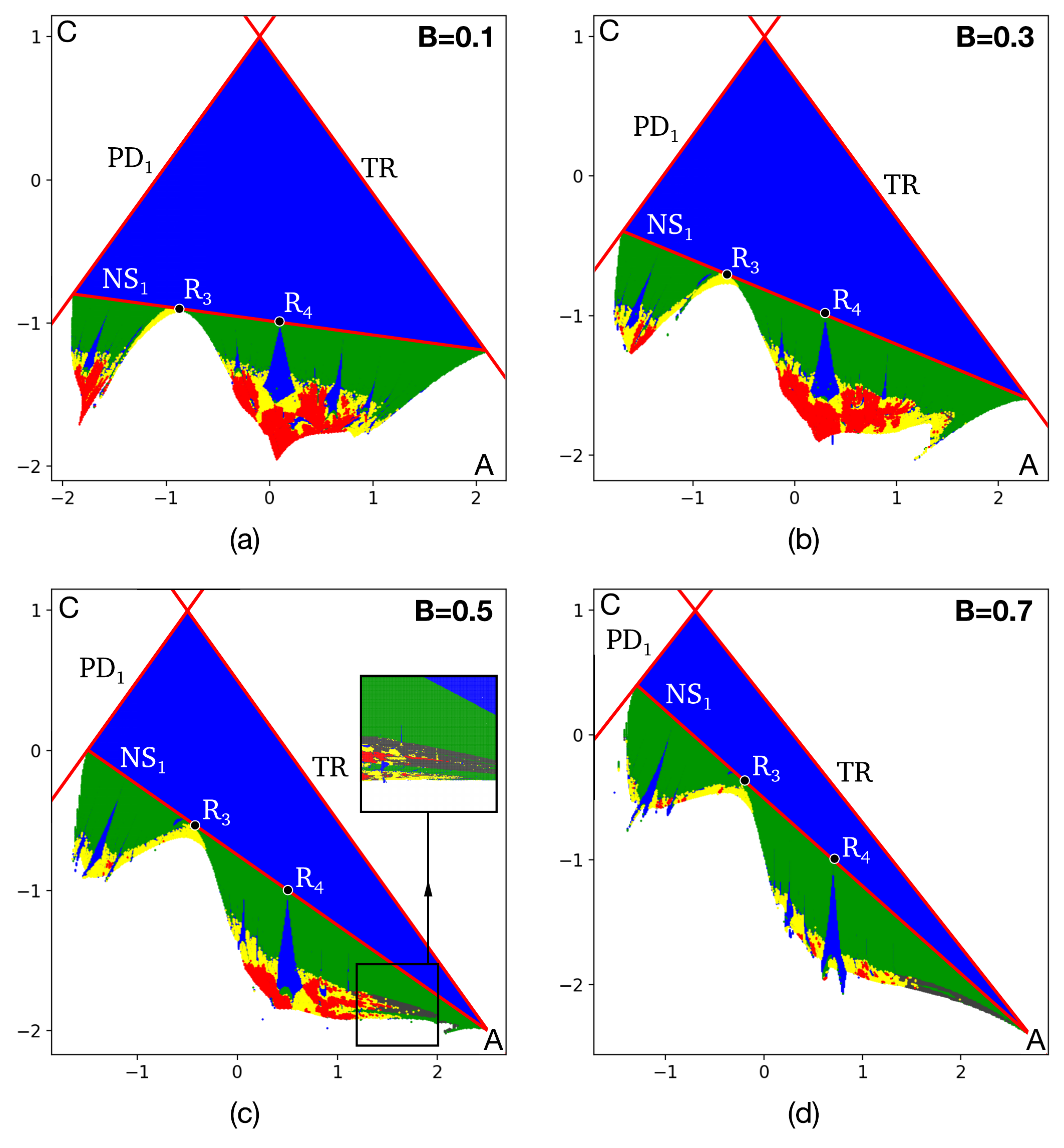}}
\end{minipage}
\caption{\footnotesize Lyapunov diagrams for map \eqref{eq:Henon3DMap} superimposed with the bifurcation curves TR, PD$_1$, and NS$_1$: (a) $B=0.1$, (b) $B=0.3$, (c) $B=0.5$, and (d) $B=0.7$. We use the following color coding: blue -- periodic regimes ($\Lambda_1 < 0$); green -- quasiperiodic regimes ($\Lambda_1 = 0, \Lambda_2 < 0$); yellow -- chaotic attractors with $\Lambda_1 > 0, \Lambda_2 < 0$; gray -- ``flow-like'' chaotic attractors, when $\Lambda_1 > 0, \Lambda_2 \approx 0$; and red -- hyperchaotic attractors with $\Lambda_1 > 0, \Lambda_2 > 0$.}
\label{ris:Fig10}
\end{figure}

In Figure~\ref{ris:Fig10} we show diagrams of Lyapunov exponents $\Lambda_1 \geq \Lambda_2 > \Lambda_3$ (Lyapunov diagrams) for various values of the Jacobian $B$ ($B \in \{0.1,0.3,0.5,0.7\}$).\footnote{On these Lyapunov diagrams we plot only regimes associated with the fixed point $O_+$. Above the curve TR, corresponding to a transcritical saddle-node bifurcation, the point $O_-$ becomes stable (this point is swapped by the stability with $O_+$). However, we do not color the corresponding part of $(A,C)$-parameter plane.} For calculation of Lyapunov exponents we take an orbit on an attractor and estimate the exponents along its $10^6$ iterations using the standard scheme \cite{BGGS80}. Depending on values of LE, we use the following color coding: blue -- for periodic regimes ($\Lambda_1 < 0$), green -- for quasiperiodic regimes ($\Lambda_1 = 0, \Lambda_2 < 0$), yellow -- for chaotic attractors ($\Lambda_1 > 0, \Lambda_2 < 0$), gray -- for ``flow-like'' chaotic attractors ($\Lambda_1 > 0, |\Lambda_2| < 0.002$) and red -- for hyperchaotic attractors ($\Lambda_1 > 0, \Lambda_2 > 0$).

\begin{figure}[h!]
\begin{minipage}[h]{1\linewidth}
\center{\includegraphics[width=1\linewidth]{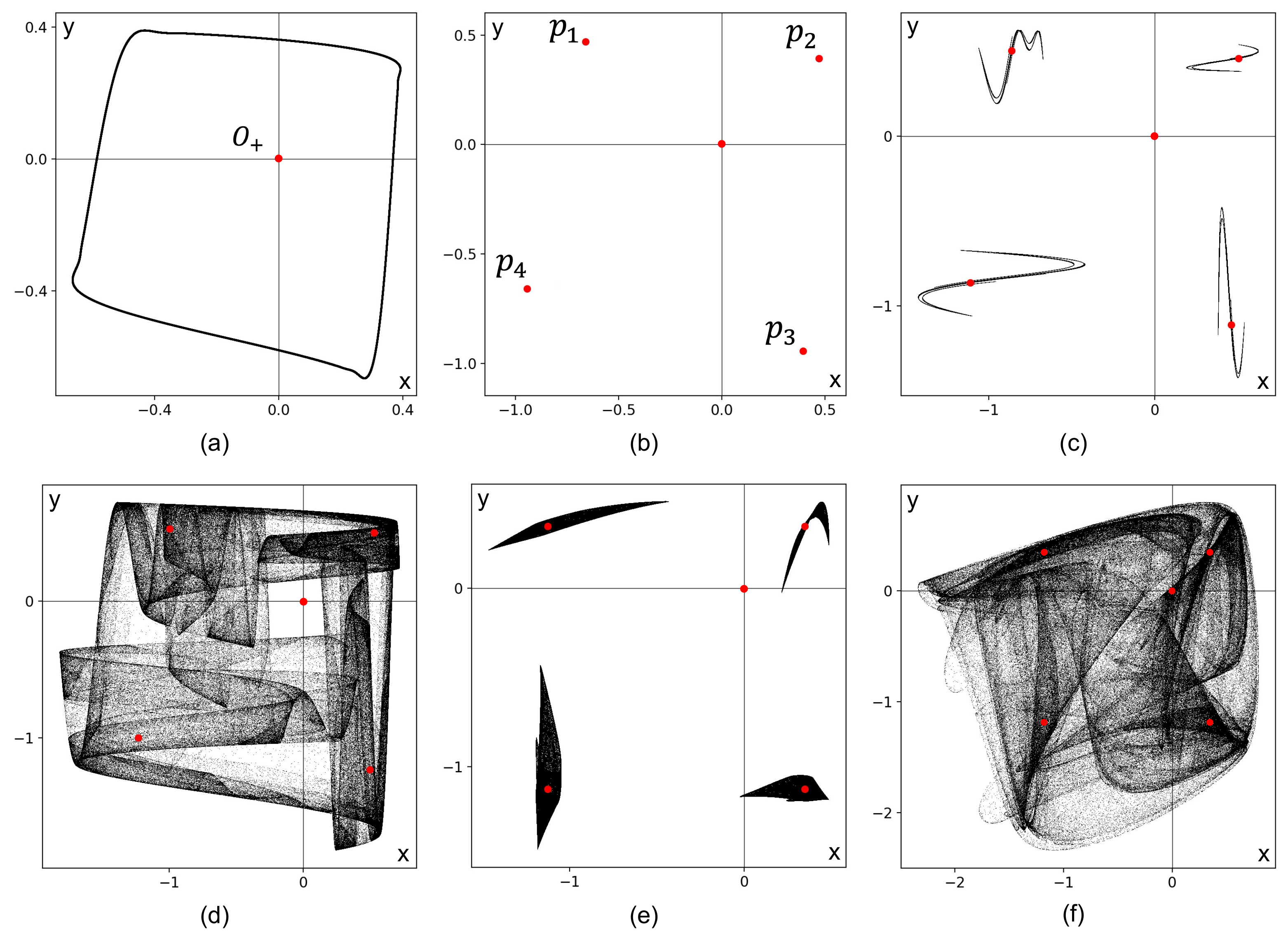}}
\end{minipage}
\caption{\footnotesize Phase portraits of various attractors of map \eqref{eq:Henon3DMap}: (a) $A=0, B=0.1, C=-1.2$; (b) $A=0, B=0.1, C=-1.4$; (c) $A=0, B=0.1, C=-1.525$; (d) $A=0, B=0.1, C=-1.62$; (e) $A=0.5, B=0.5, C=-1.745$; (f) $A=0.5, B=0.5, C=-1.78$.}
\label{ris:Fig11}
\end{figure}

Phase portraits of various attractors (for different parameters) are shown in Figure~\ref{ris:Fig11}. In Fig.~\ref{ris:Fig11}a we demonstrate the stable invariant curve which appears after the supercritical Neimark-Sacker bifurcations with the fixed point $O_+$; in Fig.~\ref{ris:Fig11}b -- the stable period-4 orbit which emerges as a resonance on this curve; in Fig.~\ref{ris:Fig11}c -- the four-component H\'enon-like attractor; in Fig.~\ref{ris:Fig11}d -- the hyperchaotic H\'enon-like attractor on the base of period-4 saddle orbit, in Fig.~\ref{ris:Fig11}e -- the hyperchaotic four-component Shilnikov attractor; and in Fig.~\ref{ris:Fig11}f -- the hyperchaotic Shilnikov attractor containing the fixed point $O_+$.

As one can see in Fig.~\ref{ris:Fig10}, hyperchaotic attractors occupy large regions under the NS$_1$-curve of the supercritical Neimark-Sacker bifurcation. After this bifurcation, the point $O_+$ becomes saddle-focus of (1,2)-type and a stable invariant curve $L$ is born in its neighborhood, see e.g. Fig.~\ref{ris:Fig11}a. Depending on the rotation number $\rho_N$, this curve can be resonant (if $\rho_N$ is rational, i.e., $\rho_N = p/q$) or ergodic (if $\rho_N$ is irrational). In the parameter plane, regions corresponding to the ergodic stable curve alternate with Arnold tongues originating from the curve of the Neimark-Sacker bifurcation where $\rho_N=p/q$. Inside Arnold tongues (close enough to the curve NS$_1$) the stable invariant curve is resonant, see e.g. Fig.~\ref{ris:Fig11}b. It is formed by the closure of the unstable invariant manifold of the period-$q$ saddle orbit emerging together with the stable one on the boundaries of the corresponding Arnold tongue.

The largest tongues correspond to the resonances with small $q$. Among them, the so-called strong resonances 1:3 and 1:4 are the most interesting and important. A tongue corresponding to the 1:3 resonance originates from a point $(A,C) = (B-1,B-1)$. It gives quite thin regions with stable dynamics, see Fig.~\ref{ris:Fig10}. However, inside it the emergence of chaotic \cite{GGS12, GGKT14} and, even, hyperchaotic attractors is possible, see more detail in Section~\ref{sec:VarShil}. On the contrary, a tongue corresponding to the 1:4 resonance gives the largest area with stable dynamics (see again Lyapunov diagrams in Fig.~\ref{ris:Fig10}). It originates from a point $(A,C) = (B,-1)$. The stable period-4 orbit inside this tongue can give rise to different types of chaotic and, even, hyperchaotic homoclinic attractors containing this orbit. If the period-4 orbit undergoes a period-doubling bifurcation, we can obtain H\'enon-like attractors, see Fig.~\ref{ris:Fig11}c, and, finally, hyperchaotic H\'enon-like attractors, see Fig.~\ref{ris:Fig11}d. If the period-4 orbit undergoes the supercritical Neimark-Sacker bifurcation, giving a four-component stable invariant curve, we, then, can observe a four-component hyperchaotic Shilnikov attractor, see Fig.~\ref{ris:Fig11}e, and, finally, hyperchaotic Shilnikov attractor containing the saddle-focus fixed point $O_+$. It is important to note that similar attractors and transitions to them are observed also in other Arnold tongues. In the framework of this paper, we study bifurcations associated only with the strong resonances 1:3 and 1:4. In particular, we show that hyperchaotic attractors on the base of the period-4 resonant orbit appear in accordance with the scenarios presented in Sec.~\ref{sec:Scenarios}.

\section{Hyperchaos in map \eqref{eq:Henon3DMap} via cascades of period-doubling bifurcations with periodic saddle orbits} \label{sec:SmallB}

In this section, we study one of possible mechanisms for the appearance of hyperchaotic attractors in map \eqref{eq:Henon3DMap} with sufficiently small values of the Jacobian $B$.

\subsection{The case $B=0$ (numerical analysis of the 2D Mir\'a map)}

\begin{figure}[h!]
\begin{minipage}[h]{1\linewidth}
\center{\includegraphics[width=1\linewidth]{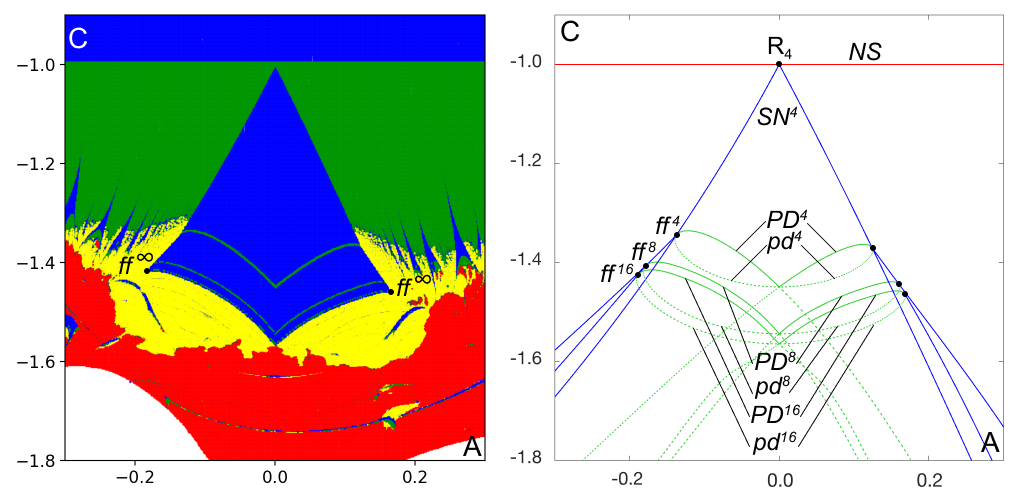}}
\end{minipage}
\caption{\footnotesize (a) Fragment of Lyapunov diagram near the 1:4 resonance for map \eqref{eq:Henon3DMap} when $B=0$; (b) corresponding bifurcation diagram on which saddle-node bifurcation curves SN$^4$, SN$^8$, and SN$^{16}$ are colored in blue and period-doubling bifurcation curves PD$^4$, PD$^8$, PD$^{16}$, pd$^4$, pd$^8$, and pd$^{16}$ -- in green (solid curves are used for bifurcations with stable orbits, while dashed -- with saddle ones); the points $ff^4$, $ff^8$, and $ff^{16}$ correspond to fold-flip bifurcations; $ff^\infty$ is a limit point for the fold-flip bifurcations.}
\label{ris:Fig12}
\end{figure}

Let us start with the case $B=0$, when this map degenerates to the two-dimensional Mir\'a endomorphism. Lyapunov diagram near the 1:4 resonance, occurring at the point $(A,C) = (0, -1)$, and the corresponding part of the bifurcation diagram obtained with help of the MatContM package \cite{MatContM13,MatContM17,MatContM20} are presented in Fig.~\ref{ris:Fig12}. One can see that hyperchaotic attractors can appear right below the region with a stable periodic dynamics. Here, we explain the organization of bifurcation curves associated with the 1:4 resonance and the scenario for the appearance of hyperchaotic attractors in this case.

We denote the stable period-4 resonant orbit $P^4 = (p_1, p_2, p_3, p_4)$, and the corresponding saddle orbit (of (2,1)-type) $S^4 = (s_1, s_2, s_3, s_4)$. The left and right boundaries of the stability region for $P^4$ are formed by a pair of saddle-node bifurcation cures SN$^4$ originating from the codimension-two point R$_4$ where $O_+$ has a pair of multipliers $e^{\pm i \pi/2}$. The bottom boundary consists of two curves PD$^4$ corresponding to the supercritical period-doubling bifurcation of $P^4$. The curve PD$^4$ touches the curve SN$^4$ at a fold-flip bifurcation point $ff^4$, where the $P^4$ has a pair of multipliers $(+1,-1)$. The same as for the fixed points of map \eqref{eq:MiraMap}, here we observe the fourth case from \cite{kuznetsov2004fold} (see its unfolding in Fig.~\ref{ris:Fig8}d). Namely, above this point, on the curves PD$^4$, the period-doubling bifurcation occurs with the stable periodic orbit $P^4$, while below this point, on the curves pd$^4$, period-doubling bifurcation occurs with the saddle period-4 orbit $S^4$. Respectively, the period-doubling bifurcation transforms the stable orbit $P^4$ to the saddle of type (2,1) and the saddle orbit $S^4$ to the saddle of type (1,2).

Very close to the point $ff^4$, there is one more important codimension-two bifurcation, significantly contributing to the organization of the bifurcation diagram. Namely, on the curve pd$^4$, we observe a point of the degenerate period-doubling bifurcation which gives rise to a curve of a saddle-node bifurcation SN$^8$ corresponding to the birth of a pair of stable and saddle period-8 orbits $P^8$ and $S^8$. The fragment of the upper boundary of the existence region for $P^8$ is formed by the curve PD$^4$. The bottom boundary of this region consists (as for the orbit $P^4$) of two pieces PD$^8$ corresponding to the supercritical period-doubling bifurcation. Similar to the case of period-4 orbit $P^4$, the period-doubling bifurcation curve PD$^8$ touches the curve SN$^8$ at the fold-flip bifurcation point $ff^8$. The bottom branch of this curve (pd$^8$) corresponds to the period-doubling bifurcation with the saddle orbit $S^8$. In its turn, close to the point $ff^8$, on the curve pd$^8$, we again observe the degenerate period-doubling bifurcation giving rise to a saddle-node bifurcation curve SN$^{16}$ and so on.

Numerical experiments show that an infinite sequence of the fold-flip points $ff^{4}$, $ff^{8}$, $ff^{16}, \dots$ accumulates to some point $ff^{\infty}$ belonging to the boundary between periodic and chaotic dynamics, see Fig.~\ref{ris:Fig12}a. This cascade of the fold-flip bifurcations gives rise to a pair of cascades of period-doubling bifurcations: the first -- with the stable periodic orbit $P^4$ (on the curves PD$^{4}$, PD$^{8}$, PD$^{16}, \dots$) and the second -- with the saddle periodic orbit $S^4$ (on the curves pd$^{4}$, pd$^{8}$, pd$^{16}, \dots$). In its turn, the periodic orbits $P^4, P^8, \dots$ which become saddles of type (2,1) after the corresponding period-doubling bifurcations, undergo again cascades of period-doubling bifurcations transforming them from the saddles of (2,1)-type to the saddles of (1,2)-type. The dashed lines PD$^{4}$, PD$^{8}$, PD$^{16}, \dots$ correspond to the first period-doubling bifurcations along these cascades.

We note that the curves PD$^{4}$ and pd$^{4}$, PD$^{8}$ and pd$^{8}$, etc., have additional common codimension-two points where the corresponding periodic orbit has the pair of multipliers $(-1,-1)$. These points belong to the line $(B=0, A=0)$. Passing through this line, both stable $P^4, P^8, \dots$ and saddle $S^4, S^8, \dots$ periodic orbits become saddles of type (1,2) simultaneously. The cascade of these codimension-two bifurcations accumulates to the so-called ``double Feigenbaum point'' \cite{KKS97}. Hyperchaotic attractor along this pathway appears exactly right below the blue colored region with periodic dynamics. This phenomenon has a simple explanation. When parameters $A$ and $B$ vanish\footnote{The evolution of chaotic dynamics with the transition from $A=0$ and $B=0$ to small values of these parameters can be rigorously traced using the anti-integrable limit approach \cite{aubry1995anti, qin2001chaotic, misiurewicz2001topological, li2006topological}. The applicability of this method for the analysis of 3D H\'enon maps was shown in \cite{juang2008chaotic, hampton2022anti}. We believe that this approach can be used to prove the existence of Smale horseshoes of type (1,2) in \eqref{eq:MiraMap} for sufficiently small parameters $A$ and $B$.} the two-dimensional Mir\'a endomorphism can be rewritten as a pair of uncoupled identical parabola maps $\bar u = 1 - c u^2, \bar v = 1-c v^2$, where hyperchaotic attractors appear right after the ``last'' period-doubling bifurcation occurring at $c \approx 1.401155$ \cite{lanford1982computer}.

Similar phenomena of the organization of bifurcation curves associated with the 1:4 resonance were previously found in several works. In \cite{Govaerts2008}, the authors studied the same unfolding of the fold-flip bifurcation with the period-4 orbit and the self-intersection of the corresponding period-doubling bifurcation curve in the duopoly model of Kopel \cite{kopel1996simple}. These results were extended in the recent book \cite{kuznetsov2019numerical} where bifurcations with period-8 and period-16 orbits were also studied in this duopoly model. In Ref.~\cite{ghaziani2012resonance} the same unfolding of the fold-flip bifurcation was studied in a predator-prey map. Our studies extend these results for explanation of the birth of hyperchaotic attractors.

\subsection{The case $B=0.1$ (3D Mir\'a map with small Jacobians)}

\begin{figure}[h!]
\begin{minipage}[h]{1\linewidth}
\center{\includegraphics[width=1\linewidth]{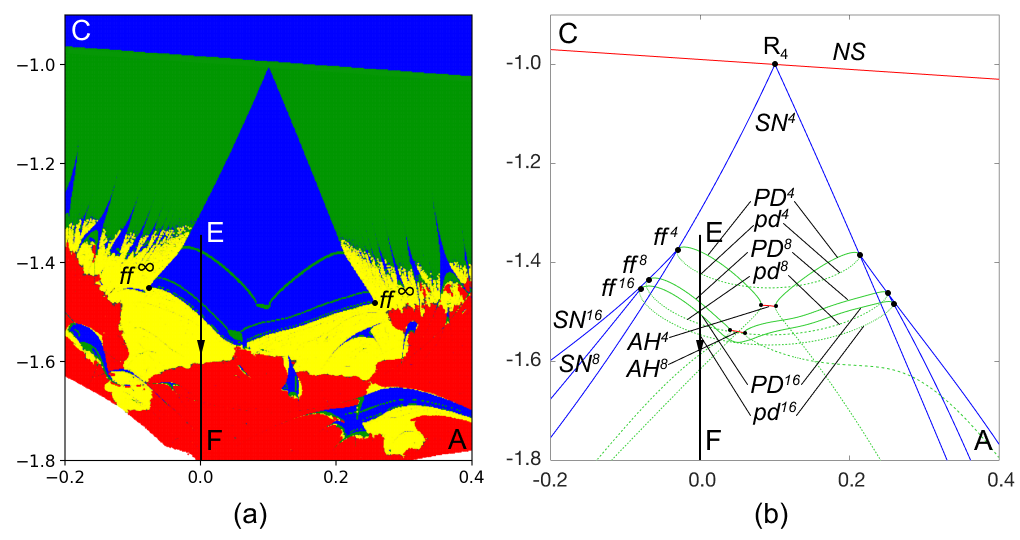}}
\end{minipage}
\caption{\footnotesize (a) Fragment of Lyapunov diagram presented in Fig.~\ref{ris:Fig10}a near the 1:4 resonance R$_4$; (b) corresponding bifurcation diagram. Here we use the same denotations for bifurcation curves and points as in Fig.~\ref{ris:Fig12}b, new red-colored curves NS$^4$ and NS$^8$ correspond to Neimark-Sacker bifurcations.}
\label{ris:Fig13}
\end{figure}

Now we fix $B=0.1$. Lyapunov diagram for this case is presented in Fig.~\ref{ris:Fig10}a. Figure~\ref{ris:Fig13} shows a zoomed region of this diagram near the 1:4 resonance, and the corresponding bifurcation diagram. Let us, first, explain changes in the organization of bifurcation curves associated with the 1:4 resonance and, then, study bifurcations along one-parameter pathways from the stable period-4 orbit $P^4$ to hyperchaotic attractors in the bottom part of this diagram.

As in the case $B=0$, the left and right boundaries of the stability region for $P^4$ are formed by a pair of saddle-node bifurcation cures SN$^4$ originating from the codimension-two point R$_4$. However the bottom boundary of this region consists here of three fragments: PD$^4$, NS$^4$, and, again, PD$^4$. The large left and right fragments correspond to the same supercritical period-doubling bifurcation, and the small middle fragment corresponds to the supercritical Neimark-Sacker bifurcation. Here, again, the curve PD$^4$ touches the curve SN$^4$ at the fold-flip point $ff^4$. We observe the same unfolding of this bifurcation as in the case $B=0$.

As for $B=0$, the stable period-8 orbit $P^8$ exists below the curve PD$^4$. The bottom boundary of this region consists (as for $P^4$) of three fragments: PD$^8$, NS$^8$, and, PD$^8$. Similar to the case of period-4 orbit $P^4$, the period-doubling bifurcation curve PD$^8$ touches the curve SN$^8$ at the fold-flip bifurcation point $ff^8$. The bottom branch of this curve (pd$^8$) corresponds to the period-doubling bifurcation with the saddle orbit $S^8$. In its turn, close to the point $ff^8$, on the curve pd$^8$, we again observe the degenerate period-doubling bifurcation giving rise to a saddle-node bifurcation curve SN$^{16}$ and so on.

One more difference between the cases $B=0$ and $B=0.1$ is that the type of period-doubling bifurcations along the curves PD$^{4}$, PD$^{8}$, $\dots$ in the second case is changed at the resonance 1:2 points where PD$^4$ intersects with NS$^4$, PD$^8$ intersects with NS$^8$, etc. As we show further, the length of curves corresponding to the Neimark-Sacker bifurcations NS$^4$, NS$^8$, $\dots$ is increased with increasing in parameter $B$.

\begin{figure}[h!]
\begin{minipage}[h]{1\linewidth}
\center{\includegraphics[width=0.8\linewidth]{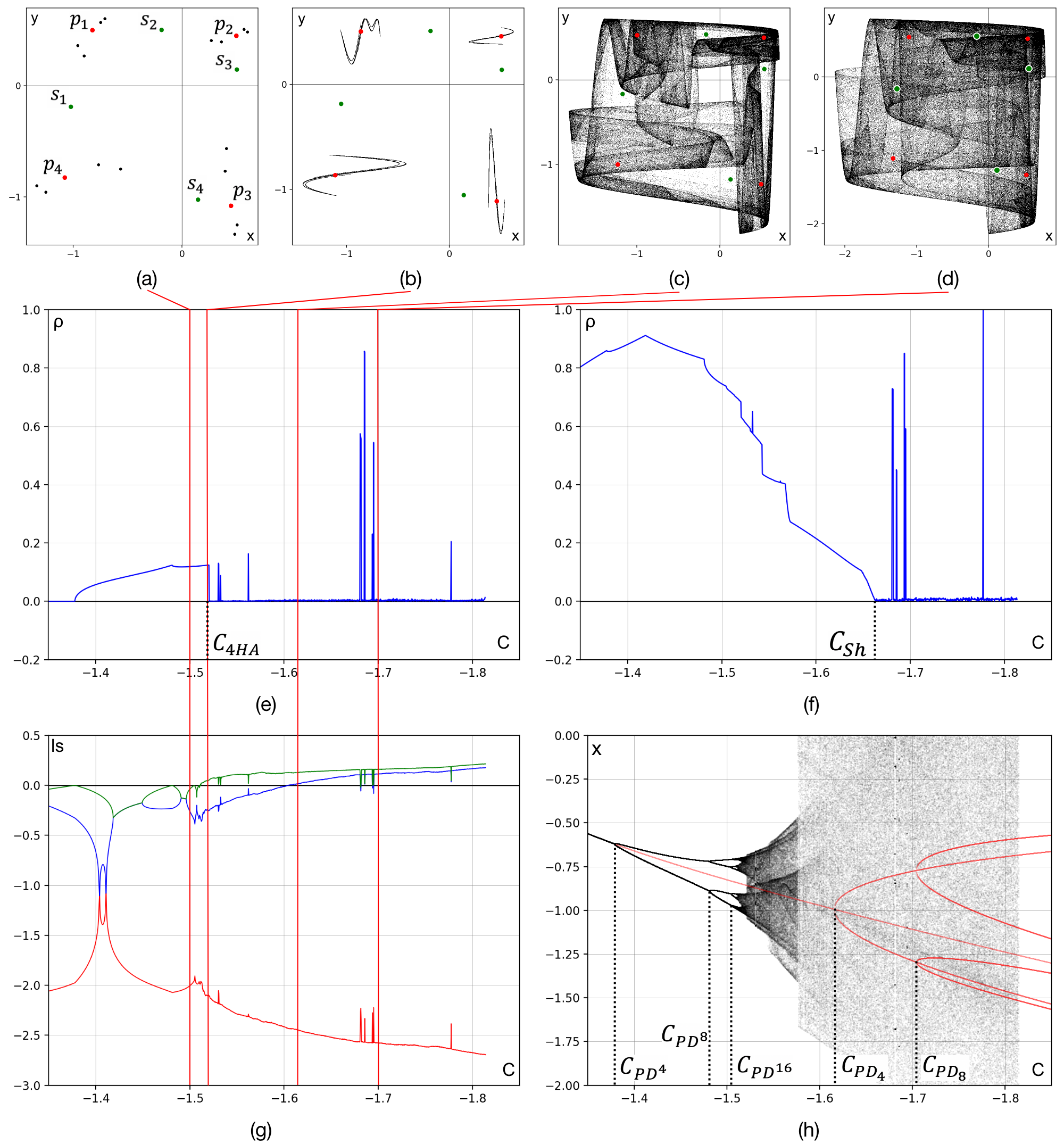}}
\end{minipage}
\caption{\footnotesize Graphs illustrating the onset of hyperchaotic attractors along the pathway EF: $(B=0.1, A=0)$. (a)--(d) phase portraits of attractors: (a) $C=-1.5$ -- period-16 stable orbit after two steps in the cascade of period-doubling bifurcation with the stable period-4 orbit $P^4: (p_1,p_2,p_3,p_4)$; (b) $C=-1.525$ -- four-component H\'enon-like attractor containing $P^4$ of (2,1)-type (LE: $\Lambda_1 = 0.073$, $\Lambda_2 = -2.156$, and $\Lambda_3 = -0.22$); (c) $C=-1.62$ -- hyperchaotic H\'enon-like attractor containing $P^4$ of (1,2)-type (LE: $\Lambda_1 = 0.132$, $\Lambda_2 = 0.027$, and $\Lambda_3 = -2.462$); (d) $C=-1.7$ -- hyperchaotic Shilnikov attractor containing the saddle-focus fixed point $O_+$ of (1,2)-type (LE: $\Lambda_1 = 0.160$, $\Lambda_2 = 0.114$, and $\Lambda_3 = -2.577$). (e) The graph of distance between the attractor and $p_1$; (f) the graph of distance between the attractor and $O_+$; (g) the graph of Lyapunov exponents $\Lambda_1, \Lambda_2$, and $\Lambda_3$ on parameter $C$; (h) bifurcation trees depicting dependency of the $x$-coordinate on parameter $C$ for every fourth iteration of the map (in black color) and continuation of $p_1$ (in red color).}
\label{ris:Fig14}
\end{figure}

\subsubsection{One-parameter bifurcation analysis.}

Along pathways transverse to the lines PD$^n$ and pd$^n$ one can expect the implementation of the scenario presented in Section~\ref{sec:period_doubling}. Let us confirm it by studying bifurcations along a one-parameter pathway EF with fixed $A = 0$. The results of corresponding analysis are shown in Figure~\ref{ris:Fig14}.

First, the stable orbit $P^4$ undergoes a cascade of period-doubling bifurcations, see phase portrait after the first two steps (slightly below the curve PD$^{16}$) in Fig.~\ref{ris:Fig14}a. The orbit $P^4$ becomes saddle of (2,1)-type after the first period-doubling bifurcation. Then, with further decrease in $C$, a cascade of heteroclinic band-merging bifurcations occurs, and, as a result, a four-component H\'enon-like attractor containing $P^4$ appears, see Fig.~\ref{ris:Fig14}b. The corresponding bifurcation tree is presented in Fig.~\ref{ris:Fig14}h. In Fig.~\ref{ris:Fig14}e depicting the graph of distance between the attractor and $p_1$ component of $P^4$, one can see that the saddle orbit $P^4$ starts to belong to the attractor (which means the H\'enon-like attractor occurrence) at $C = C_{4HA} \approx -1.52$.

The next step in the framework of hyperchaotic attractors development is a period-doubling bifurcation with the saddle orbit $P^4$ occurring on the curve $PD_4$ (at $C = C_{PD_4} \approx -1.616$). After this bifurcation, $P^4$ becomes saddle of (1,2)-type and a period-8 saddle orbit of (2,1)-type appears in its neighborhood, see the red-colored continuation tree for the $p_1$ component of $P^4$ in Fig.~\ref{ris:Fig14}h. As it can be seen in Figs.~\ref{ris:Fig14}e and Fig.~\ref{ris:Fig14}h, at $C < C_{PD_4}$ the point $p_1$ also belongs to the attractor. Moreover, as it is shown in the graph of Lyapunov exponents presented in Fig.~\ref{ris:Fig14}g, the attractor becomes hyperchaotic near the mentioned above period-doubling bifurcation. The resulting hyperchaotic attractor at $C = -1.62$ is presented in Fig.~\ref{ris:Fig14}c.

\begin{figure}[h!]
\begin{minipage}[h]{1\linewidth}
\center{\includegraphics[width=1\linewidth]{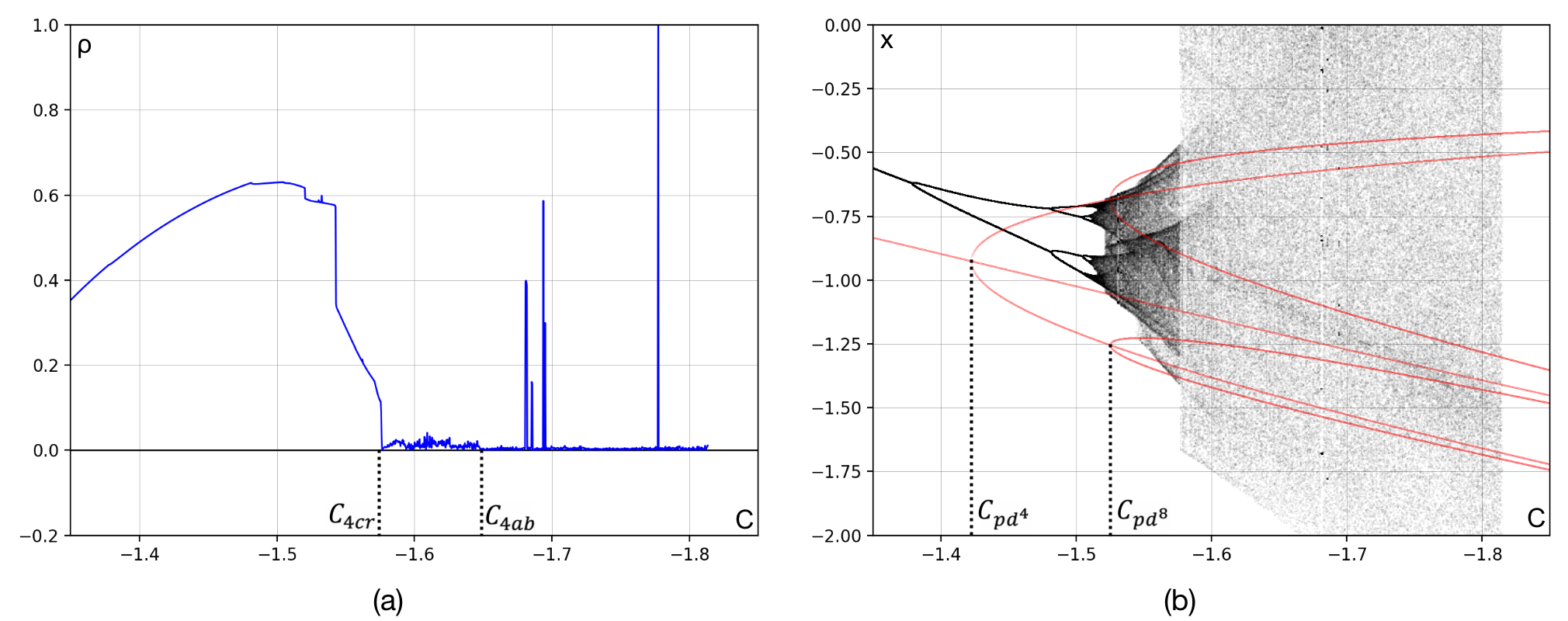}}
\end{minipage}
\caption{\footnotesize Diagrams illustrating bifurcations of the period-4 saddle orbit $S^4 = (s_1, s_2, s_3, s_4)$ which appears together with $P^4$ via a saddle-node bifurcation on the curve SN$^4$: (a) the graph of distance between the attractor and $s_1$; (b) bifurcation tree depicting dependency of $x$-coordinate on parameter $C$ for every fourth iteration of the map superimposed with the continuation tree for $s_1$.}
\label{ris:Fig15}
\end{figure}

Note that before the period-doubling bifurcation (at $C>C_{PD_4}$) with the saddle orbit $P^4$, the four-component H\'enon-like attractor collides into the one-component attractor. It happens due to the boundary crisis: the unstable invariant manifold $W^u(P^4)$, forming the attractor, begins to intersect with the two-dimensional stable invariant manifold of the period-16 saddle orbit forming the boundary of its absorbing domain. (The corresponding heteroclinic bifurcation occurs at $C = C_{4cr} \approx -1.58$, see the jump of distance between the attractor and $s_1$ in Fig.~\ref{ris:Fig15}a). This period-16 saddle orbit appears after two period-doubling bifurcations (happened on the curves pd$^4$ and pd$^8$, see Fig.~\ref{ris:Fig13}b) with the saddle orbit $S^4$ which, as was mentioned above, goes through the full cascade of period-doubling bifurcations (see also the red-colored tree in Fig.~\ref{ris:Fig15}b). With a further decrease in $C$, at $C = C_{4ab} \approx -1.65$, the attractor absorbs the saddle orbit $S^4$ of type (1,2), see again Fig.~\ref{ris:Fig15}a.

It is worth noting that the hyperchaotic attractor presented in Fig.~\ref{ris:Fig14}c contains the orbits $P^4$ and $S^4$ but does not contain the fixed point $O_+$ which becomes a saddle-focus of (1,2)-type at $C < -0.99$. With a further decrease in $C$, the distance between the attractor and this fixed point decreases and, finally, at $C = C_{Sh} \approx -1.66$ it vanishes, see Fig.~\ref{ris:Fig14}f. As a result, a hyperchaotic Shilnikov attractor appears, see Fig.~\ref{ris:Fig14}d.

\section{Hyperchaos in map \eqref{eq:Henon3DMap} via onset of multicomponent Shilnikov attractors} \label{sec:LargeB}

In this section we study scenarios for the appearance of hyperchaotic attractors in map \eqref{eq:Henon3DMap} with not very small values of the Jacobian $B$.

\subsection{The case $B=0.5$}

Lyapunov diagram for $B=0.5$ is presented in Fig.~\ref{ris:Fig10}c. As in the previous case, hyperchaotic attractors occurs here after the destruction of the stable invariant curve $L$. However another mechanisms of the destruction of $L$ is more typical in this case.

\begin{figure}[h!]
\begin{minipage}[h]{1\linewidth}
\center{\includegraphics[width=1\linewidth]{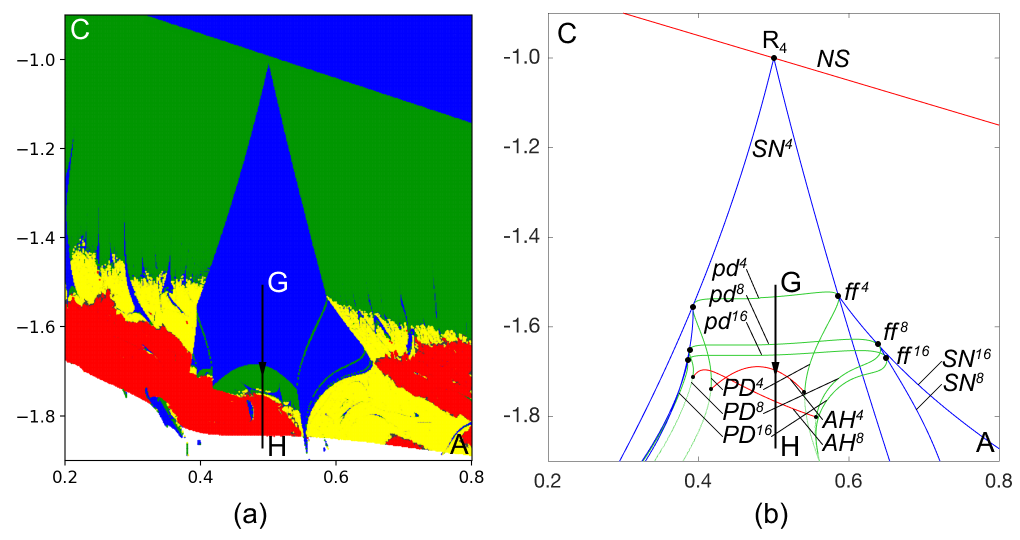}}
\end{minipage}
\caption{\footnotesize (a) $B=0.5$. Fragment of Lyapunov diagram near the 1:4 resonance R$_4$, and (b) corresponding bifurcation diagram on which we use the same denotations for bifurcation curves and points as in Fig.~\ref{ris:Fig12}b.}
\label{ris:Fig16}
\end{figure}

Figure~\ref{ris:Fig16} shows an enlarged fragment of Lyapunov diagram near the Arnold tongue for the 1:4 resonance and the corresponding bifurcation diagram where we use the same denotations for periodic orbits and bifurcation curves associated with this resonance as in Sec.~\ref{sec:SmallB}.

As in the previous case, here we also observe the cascade of fold-flip bifurcations ($ff^4$, $ff^8, \dots$) which gives rise to the pair of cascades of period-doubling bifurcations with the stable period-4 orbit $P^4$ (the curves PD$^4$, PD$^8, \dots$) and with the saddle period-4 orbit $S^4$ (the curves pd$^4$, pd$^8 \dots$). The degenerate period-doubling bifurcations occurring on the curves PD$^4$, PD$^8, \dots$ near the corresponding fold-flip points generate a cascade of saddle-node bifurcation curves SN$^4$, SN$^8, \dots$. A pair of period-($4 \cdot 2^n$) stable and saddle orbits appears at the intersection of each such curve.

However, unlike the case of small values of the Jacobian $B$, the curves pd$^4$, pd$^8 \dots$ are located above the curves PD$^4$, PD$^8, \dots$, compare Figures~\ref{ris:Fig13}b and \ref{ris:Fig16}b. More detailed study of transformation of bifurcation diagrams from the case of small values of $B$ (in particular, starting with $B=0$) to the case of not small values of $B$ seems very interesting problem for future research.

As in the case of small $B$, the region of stability for the period-4 orbit $P^4$ is bounded from below by the period-doubling bifurcation curve PD$^4$, Neimark-Sacker bifurcation curve NS$^4$, and again period-doubling bifurcation curve PD$^4$. However, the fragment NS$^4$ is much wider here comparing with the case of small $B$. Moreover, the curves NS$^4$ and PD$^4$ are organized in such a way, that typical pathways from $P^4$ to hyperchaos pass through one of the Neimark-Sacker bifurcation curves, see Fig.~\ref{ris:Fig16}.

\subsubsection{One-parameter bifurcation analysis.}

\begin{figure}[h!]
\begin{minipage}[h]{1\linewidth}
\center{\includegraphics[width=0.8\linewidth]{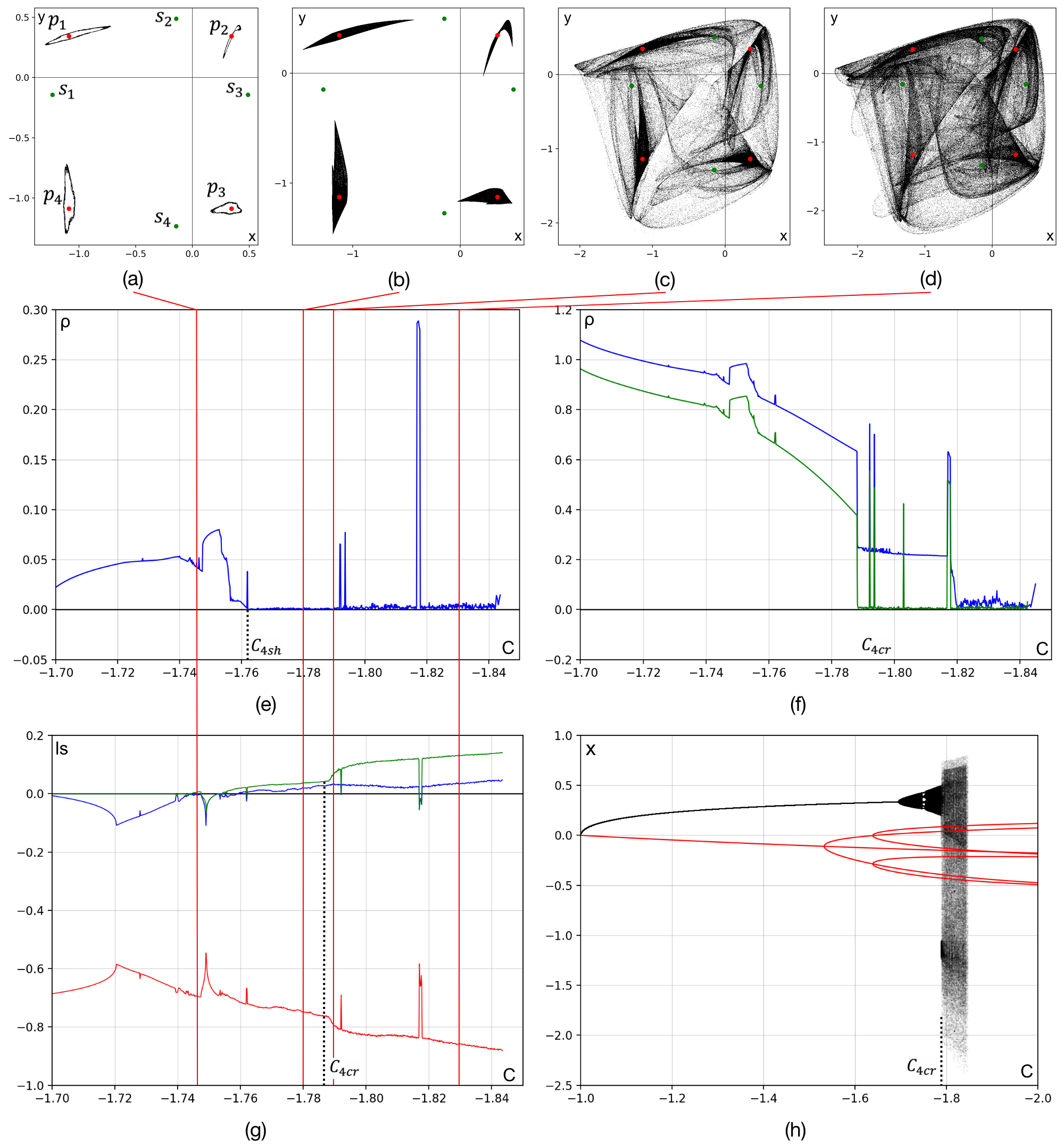}}
\end{minipage}
\caption{\footnotesize Graphs illustrating the onset of hyperchaotic attractors along the pathway GH: $(B=0.5, A=0.5)$. (a)--(d) phase portraits of attractors: (a) $C=-1.745$ -- four-component torus-chaos attractor; (b) $C=-1.78$ -- four-component hyperchaotic Shilnikov attractor containing the saddle-focus period-4 orbit $P^4$ (LE: $\Lambda_1 = 0.036$, $\Lambda_2 = 0.0187$, and $\Lambda_3 = -0.748$); (c) $C=-1.79$ -- four-component attractor transforms to the one-component attractor (LE: $\Lambda_1 = 0.072$, $\Lambda_2 = 0.0315$, and $\Lambda_3 = -0.797$); (d) $C=-1.83$ -- hyperchaotic Shilnikov attractor containing the saddle-focus fixed point $O_+$ (LE: $\Lambda_1 = 0.131$, $\Lambda_2 = 0.0346$, and $\Lambda_3 = -0.858$). (e) the graph of distance between the attractor and the $p_1$-component of $P^4$. (f) the graphs of distance between the attractor and $O_+$ (in blue color), and the attractor and the $s_1$-component of $S^4$ (in green color). (g) the graph of LE on parameter $C$. (h) bifurcation trees depicting dependency of $x$-coordinate on parameter $C$ for every fourth iteration of the map (in black color) and continuation of $s_1$ (in red color).}
\label{ris:Fig17}
\end{figure}

Further, let us fix $A=0$ (together with $B=0.5$) and study bifurcations along a pathway GH from the stable period-4 orbit $P^4$ to hyperchaotic attractors. The results of corresponding one-parameter analysis are shown in Figure~\ref{ris:Fig17}.

In contrast to the cases of small $B$ (see Sec.~\ref{sec:SmallB}), the stable period-4 orbit $P^4$ undergoes here the supercritical Neimark-Sacker bifurcation. As a result, a stable four-component invariant curve is born in the neighborhood of $P^4$ while this orbit becomes a saddle-focus of (1,2)-type. With decreasing $C$, this four-component curve breaks down and a four-component ``torus-chaos'' attractor (with only one positive Lyapunov exponent, see the graph of Lyapunov exponents in Fig.~\ref{ris:Fig17}g) appears, see Fig.~\ref{ris:Fig17}a.

\begin{figure}[!ht]
\begin{minipage}[h]{1\linewidth}
\center{\includegraphics[width=0.8\linewidth]{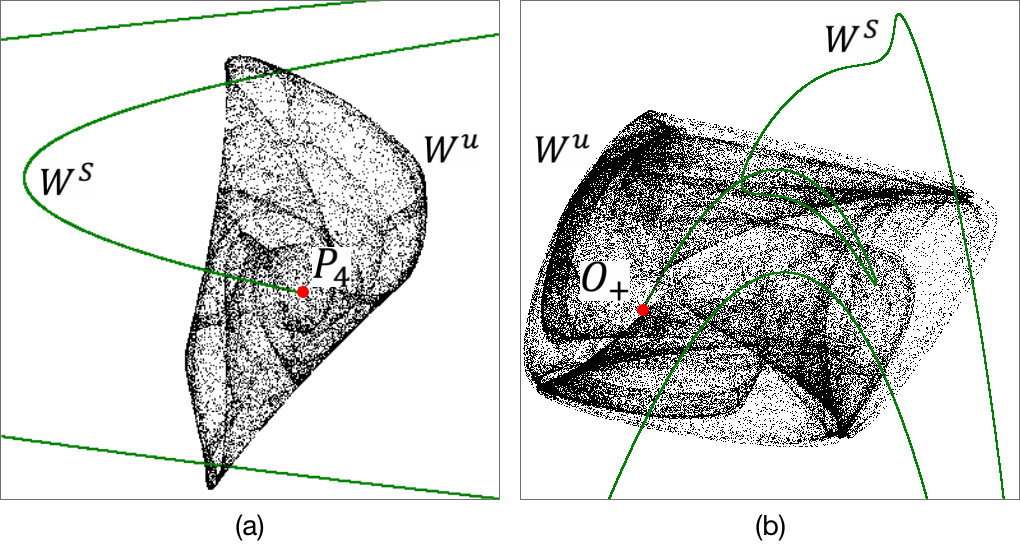}}
\end{minipage}
\caption{{\footnotesize Phase portraits of (a) one component of the hyperchaotic four-component Shilnikov attractor presented in Fig.~\ref{ris:Fig17}c, and (b) hyperchaotic Shilnikov attractor containing $O_+$ shown in Fig.~\ref{ris:Fig17}d. Green-colored curve in both panels is that branch of the stable invariant manifold $W^s$ which intersects with the unstable manifold $W^u$ forming the attractor.}}
\label{ris:Fig18}
\end{figure}

Then, according to the scenario described in \cite{GSKK19, sataev2021cascade} (see also Sec.~\ref{sec:ShilAttrScen}), this torus-chaos attractor absorbs the saddle-focus orbit $P^4$ and, as a result, a four-component hyperchaotic Shilnikov attractor appears, Fig.~\ref{ris:Fig17}b. This attractor is homoclinic, it contains $P^4$, its unstable invariant manifold $W^u$, and homoclinic points belonging to the intersection $W^u(P^4) \cap W^s(P^4)$. As it can be seen from Fig.~\ref{ris:Fig17}e, depicting the distance between the attractor and a $p_4$-component of $P^4$, the attractor starts to contain $P^4$ at $C=C_{4Sh} \approx -1.762$. Moreover, $P^4$ belongs to the attractor on a quite large interval of parameter $C$. Fig.~\ref{ris:Fig18}a depicting the phase portrait of the attractor and the stable one-dimensional manifold $W^s(P^4)$ shows that $W^u(P^4)$ and $W^s(P^4)$ intersect transversally on this interval. By Smale and Shilnikov \cite{Sm67,Sh67}, the emergence of such an intersection implies a countable many saddle-focus orbits with the two-dimensional unstable manifold inside the attractor.

The onset of the four-component Shilnikov attractor is not the final step in the framework of the development of hyperchaotic dynamics along the pathway GH. With a further decrease in $C$ (at $C = C_{4cr} \approx -1.789$), this attractor collides into the one-component attractor, see Fig.~\ref{ris:Fig17}c. As in the case described in Sec.~\ref{sec:SmallB}, this collision happens due to the boundary crisis, see the explosive growth of the attractor size in the bifurcation tree presented in Fig.~\ref{ris:Fig17}h and the jump of distance between the attractor and $s_1$-component of $S^4$ in Fig.~\ref{ris:Fig17}f (green-colored graph). However, in this case, the period-4 saddle orbit $S^4$ undergoes before the crisis of the attractor a complete cascade of period-doubling bifurcations (on the curves pd$^4$, pd$^4, \dots$ in Fig.~\ref{ris:Fig16}b), see a part of the continuation tree for the $s_1$-component of $S^4$ in Fig.~\ref{ris:Fig17}h (red-colored graph). As a result, a non-attractive hyperbolic set with the two-dimensional unstable manifolds appears on the base of this orbit. At $C = C_{4cr}$, the four-component attractor touches the stable manifolds of this set and, at $C < C_{4cr}$, merges with it into the one-component attractor. Note that this collision is also clearly visible in the graph of Lyapunov exponents presented in Fig.~\ref{ris:Fig17}g, where one can observe the jump of the second Lyapunov exponent (the so-called jump of hyperchaoticity \cite{StanKazGon20}) at $C = C_{4cr}$ associated with this collision.

The one-component hyperchaotic attractor presented in Fig.~\ref{ris:Fig17}c contains both periodic saddle orbits $P^4$ and $S^4$ which have the two-dimensional unstable manifold but does not contain the saddle-focus fixed point $O_+$. With a further decrease in $C$, the distance between the attractor and this point decreases and finally, at $C = C_{ab} \approx -1.820$, it vanishes, see the blue-colored graph in Fig.~\ref{ris:Fig17}f. As a result, a hyperchaotic Shilnikov attractor containing the point $O_+$ appears, see Fig.~\ref{ris:Fig17}d. It contains infinitely many saddle-focus, as well as saddle periodic orbits with two-dimensional unstable manifolds. Figure~\ref{ris:Fig18}b shows the transversal homoclinic structure for this attractor, additionally confirming the inclusion of $O_+$ to it.

Further decrease in $C$ leads to the destruction of the homoclinic Shilnikov attractor which happens via a boundary crisis: the unstable manifold $W^u(O_+)$ forming the attractor intersects with the stable two-dimensional manifold $W^s(O_-)$ bounding its absorbing domain.

We would like to note that a similar transition to hyperchaos is observed also along many other pathways from the stable fixed point $O_+$ to hyperchaotic attractors. Passing through other Arnold tongues, corresponding e.g. to a $p/q$-resonance, one can observe a transition from a stable period-$q$ orbit to a Shilnikov attractor containing this orbit which becomes a saddle-focus of (1,2)-type after the corresponding Neimark-Sacker bifurcation. Then, the $q$-component Shilnikov attractor collides into a one-component attractor which, with a further decrease in $C$, can absorb the saddle-focus fixed point $O_+$ before the crisis.

\begin{figure}[h!]
\begin{minipage}[h]{1\linewidth}
\center{\includegraphics[width=1\linewidth]{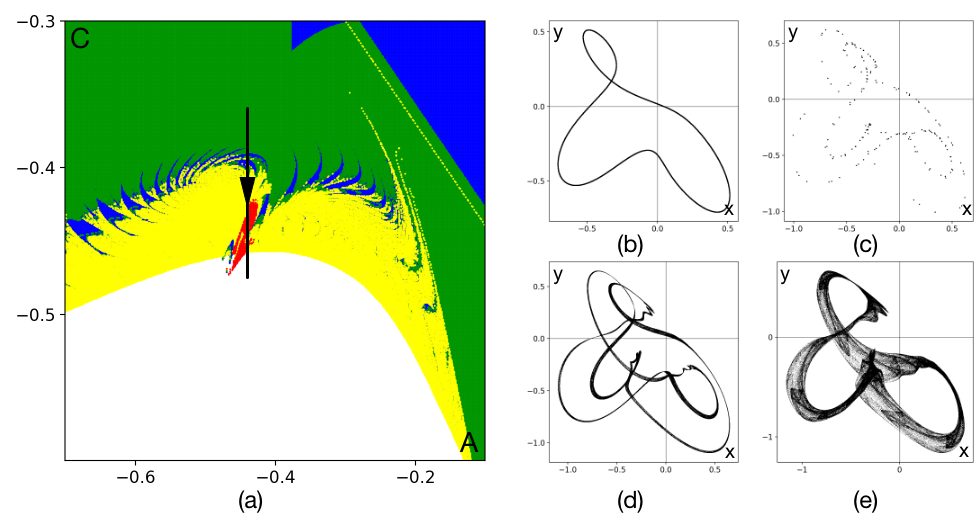}}
\end{minipage}
\caption{\footnotesize (a) Fragment of Lyapunov diagram near the 1:3 resonance, $B = 0.7$. (b)--(e) Phase portraits of attractors along the pathway with fixed $A=-0.44$: (b) $C=-0.36$, the stable invariant curve $L$ deforms near the 1:3 resonance; (c) long-period resonance appears inside the corresponding small Arnold tongue; (d) $C=-0.41$, torus-chaos attractor (LE: $\Lambda_1 = 0.017$, $\Lambda_2 = -0.011$, and $\Lambda_3 = -0.363$) (e) $C=-0.44$, hyperchaotic attractor (LE: $\Lambda_1 = 0.015$, $\Lambda_2 = 0.004$, and $\Lambda_3 = -0.376$).}
\label{ris:Fig19}
\end{figure}

However, it is not the case for transition near the strong 1:3 resonance. We do not observe a stable period-3 orbit inside the corresponding tongue despite the possible existence (e.g. for $B=0.7$) of hyperchaotic attractors below it, see Fig.~\ref{ris:Fig10}d. The enlarged fragment of the corresponding Lyapunov diagram is presented in Figure~\ref{ris:Fig19}a. Let us briefly describe main stages of the development of hyperchaotic attractors along a vertical pathway passing through $A=-0.44$. The stable invariant curve $L$ deforms near the 1:3 resonance, see Fig.~\ref{ris:Fig19}b. Then, a higher periodic resonance occurs inside the corresponding Arnold tongue, see Fig.~\ref{ris:Fig19}c. With a further decrease in $C$, a multi-component invariant curve, appears from this periodic orbit under the supercritical Neimark-Sacker bifurcation. Soon, this curve breaks down giving torus-chaos attractor, see Fig.~\ref{ris:Fig19}d. Finally, this strange attractor becomes hyperchaotic, see Fig.~\ref{ris:Fig19}e.

\section{Chaotic attractors with zero second Lyapunov exponent} \label{sec:VarShil}

In Sections~\ref{sec:SmallB} and~\ref{sec:LargeB} we have studied the appearance of homoclinic attractors in map \eqref{eq:Henon3DMap} for both cases of small and not very small values of the Jacobian $B$. In these cases the evolution of attractors is terminated via the boundary crisis. We have shown that, before this, the attractor can absorb the fixed point $O_+$ of the saddle-focus (1,2)-type, which gives the formal possibility for the attractor to become hyperchaotic.

However, chaotic attractors do not always result in hyperchaos in the map under consideration. In this section, we show that for sufficiently large values of the Jacobian $|B| < 1$ chaotic attractors (including Shilnikov ones) can be flow-like, i.e., their second Lyapunov exponent can be indistinguishable from zero in numerical experiments in quite large regions of the parameter space. This in itself is not a new phenomenon. Flow-like chaotic attractors were previously observed in the three-dimensional H\'enon map \cite{GOST05}, in the Lorenz-84 model \cite{BSV05}, in some class of 3D diffeomorphisms of the torus \cite{BSV10}, in nonholonomic models of Celtic stone \cite{GGKS19} and Chaplygin top \cite{BKS16}, in models of identical globally coupled oscillators \cite{grines2022origin} and in other systems. Also we give an explanation for this phenomenon.

\begin{figure}[h!]
\begin{minipage}[h]{1\linewidth}
\center{\includegraphics[width=1\linewidth]{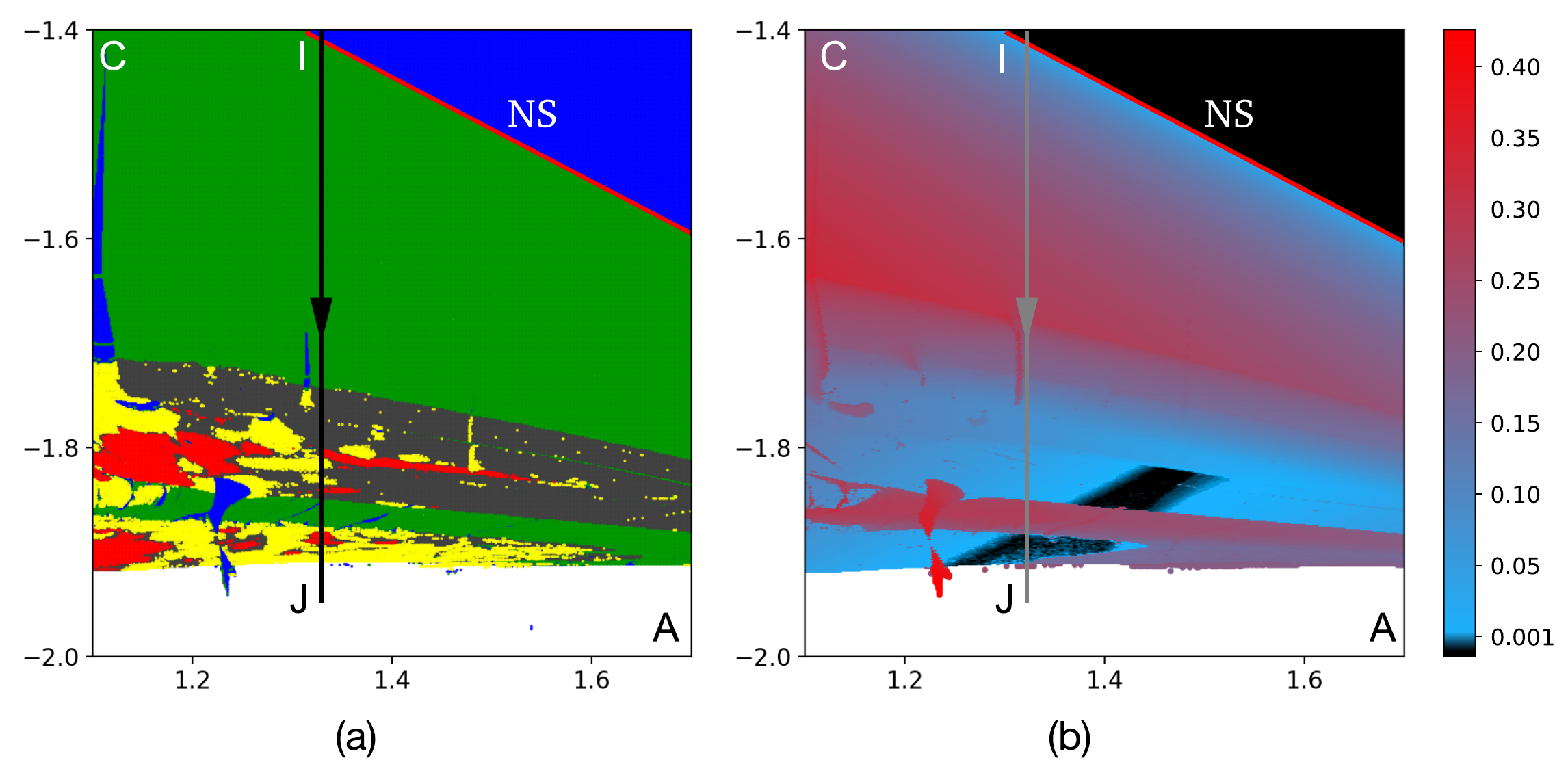}}
\end{minipage}
\caption{\footnotesize ($A,C$)-parameter diagrams for map \eqref{eq:Henon3DMap}, $B=0.5$. (a) Lyapunov diagram. (b) Diagram of the distance between the attractor and the fixed point $O_+$, in black colored regions this distance is less than 0.001. In the top-right region the point $O_+$ is asymptotically stable; a homoclinic attractors containing the saddle-focus point $O_+$ of type (1,2) appear in the bottom-center region, colored in black in Fig. (b), it can be either hyperchaotic or not.}
\label{ris:Fig20}
\end{figure}

Figure~\ref{ris:Fig20} shows two diagrams for map \eqref{eq:Henon3DMap} with $B=0.5$ on the ($A,C$)-parameter plane. Panel (a) is the enlarged fragment of the right-bottom part of the Lyapunov diagram presented in Fig.~\ref{ris:Fig10}c, and panel (b) is the corresponding piece of a distance diagram depicting the distance between the attractor and the fixed point $O_+$. For the Lyapunov diagram we use the same color coding as in Sec.~\ref{sec:MiraMapBiff}. Black color in the distance diagram corresponds to a small distance between the attractor and $O_+$ (the minimal distance less than 0.001 after $10^6$ iterations of a point taken on the attractor), see the legend to the right of panel (b). The top-right corner of Figs.~\ref{ris:Fig20} corresponds to the stability region of the point $O_+$, homoclinic attractors containing the saddle-focus fixed point $O_+$ of type (1,2) appear in the bottom-center region, colored in black at panel (b). As one can see in the Lyapunov diagram, strange attractors here can be of three possible types: strongly dissipative (yellow color, $\Lambda_1 > 0$, $\Lambda_2 < 0$), hyperchaotic (red color, $\Lambda_1 > 0$, $\Lambda_2 > 0$), and flow-like (gray color, $\Lambda_1 > 0$, $\Lambda_2 \approx 0$).

\begin{figure}[h!]
\begin{minipage}[h]{1\linewidth}
\center{\includegraphics[width=1\linewidth]{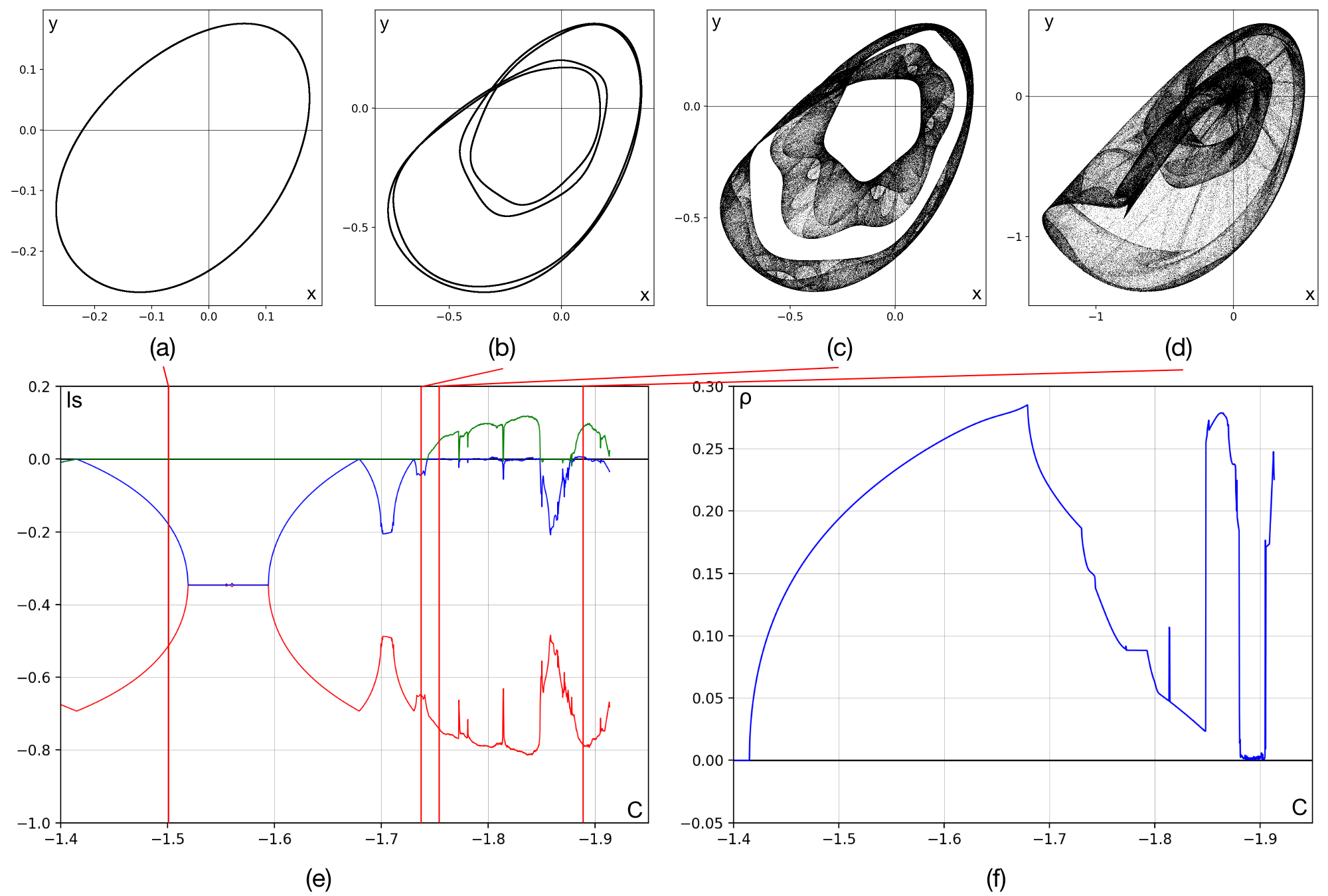}}
\end{minipage}
\caption{\footnotesize The evolution of attractors along the pathway IJ: $(B=0.5, A=1.33)$. (a)--(d) phase portraits of attractors: (a) $C=-1.5$ -- the stable invariant curve $L$; (b) $C=-1.732$ -- four-round invariant curve; (c) $C=-1.75$ -- torus-chaos attractor (LE: $\Lambda_1 = 0.033$, $\Lambda_2 = -3.65 \cdot 10^{-5}$, and $\Lambda_3 = -0.726$); (d)~$C=-1.89$ -- ``weakly'' hyperchaotic Shilnikov attractor (LE: $\Lambda_1 = 0.09$, $\Lambda_2 = 0.005$, and $\Lambda_3 = -0.788$). (e)~the graph of Lyapunov exponents $\Lambda_1, \Lambda_2$, and $\Lambda_3$ on parameter $C$. (f) the graph of the distance between the attractor and the fixed point $O_+$.}
\label{ris:Fig21}
\end{figure}

Let us further fix $A=1.33$ and study the scenario of strange attractors appearance along the pathway IJ: $(A = 1.33, B=0.5)$. Results of the corresponding bifurcation analysis are shown in Figure~\ref{ris:Fig21}. At the beginning, the stable fixed point $O_+$ is the only attractor of the map. Then, on the curve NS, this point undergoes the supercritical Neimark-Sacker bifurcation after which $O_+$ becomes saddle-focus of type (1,2) and the stable invariant curve $L$ appears in its neighborhood, see Fig.~\ref{ris:Fig21}a. Then, this curve undergoes two period-doubling (length-doubling) bifurcations. After the first period-doubling this curve becomes saddle and a stable doubled (2-round) invariant curve $L^2$ appears in its neighborhood. In its turn, curve $L^2$ also undergoes a period-doubling bifurcation, it becomes saddle and a four-round stable invariant curve appears, see Fig.~\ref{ris:Fig21}b. With further decrease in $C$, the four-round invariant curve breaks down, and, as a result, a torus-chaos attractor, possessing near-zero second Lyapunov exponent, is born, see Fig.~\ref{ris:Fig21}c.

\begin{figure}[h!]
\begin{minipage}[h]{1\linewidth}
\center{\includegraphics[width=1\linewidth]{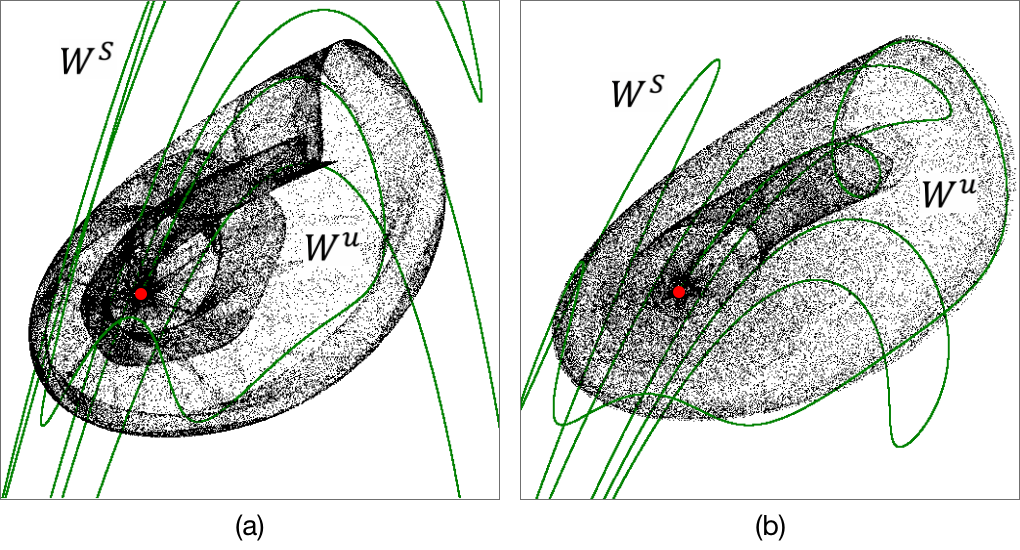}}
\end{minipage}
\caption{\footnotesize Phase portraits for different discrete Shilnikov attractors containing $O_+$: (a) $B=0.5, A=1.33, C=-1.89$ (LE: $\Lambda_1 = 0.09$, $\Lambda_2 = 0.005$, and $\Lambda_3 = -0.788$), (b) $B=0.7, A=1.55, C=-2.03$ (LE: $\Lambda_1 = 0.107$, $\Lambda_2 \approx 0$, and $\Lambda_3 = -0.463$). Green-colored curve is that branch of stable invariant manifold of $O_+$ which forms the homoclinic intersection with the unstable manifold of $O_+$.}
\label{ris:Fig22}
\end{figure}

The graph of distance between the attractor and the fixed point $O_+$, presented in Fig.~\ref{ris:Fig21}f, shows that the attractor starts to contain $O_+$ at $C \approx -1.882$. Figure~\ref{ris:Fig22}a, showing the transverse homoclinic structure for $O_+$, additionally confirms the inclusion of $O_+$ to the attractor. As it can be seen from the graph of Lyapunov exponents (Fig.~\ref{ris:Fig21}e), the attractor becomes hyperchaotic on the interval $C \in (-1.89,-1.88)$. However, since the second Lyapunov exponent is slightly positive ($0 < \Lambda_2 < 0.006)$, hyperchaos here is very weak.

Moreover, the second Lyapunov exponent vanishes for chaotic attractors existing at some sufficiently large, open regions of the parameter space adjacent to the regions with homoclinic attractors. In order to see it, just superimpose the black-colored region in the bottom-center part of Fig.~\ref{ris:Fig20}b with the corresponding part of the Lyapunov diagram.

\begin{figure}[h!]
\begin{minipage}[h]{1\linewidth}
\center{\includegraphics[width=1\linewidth]{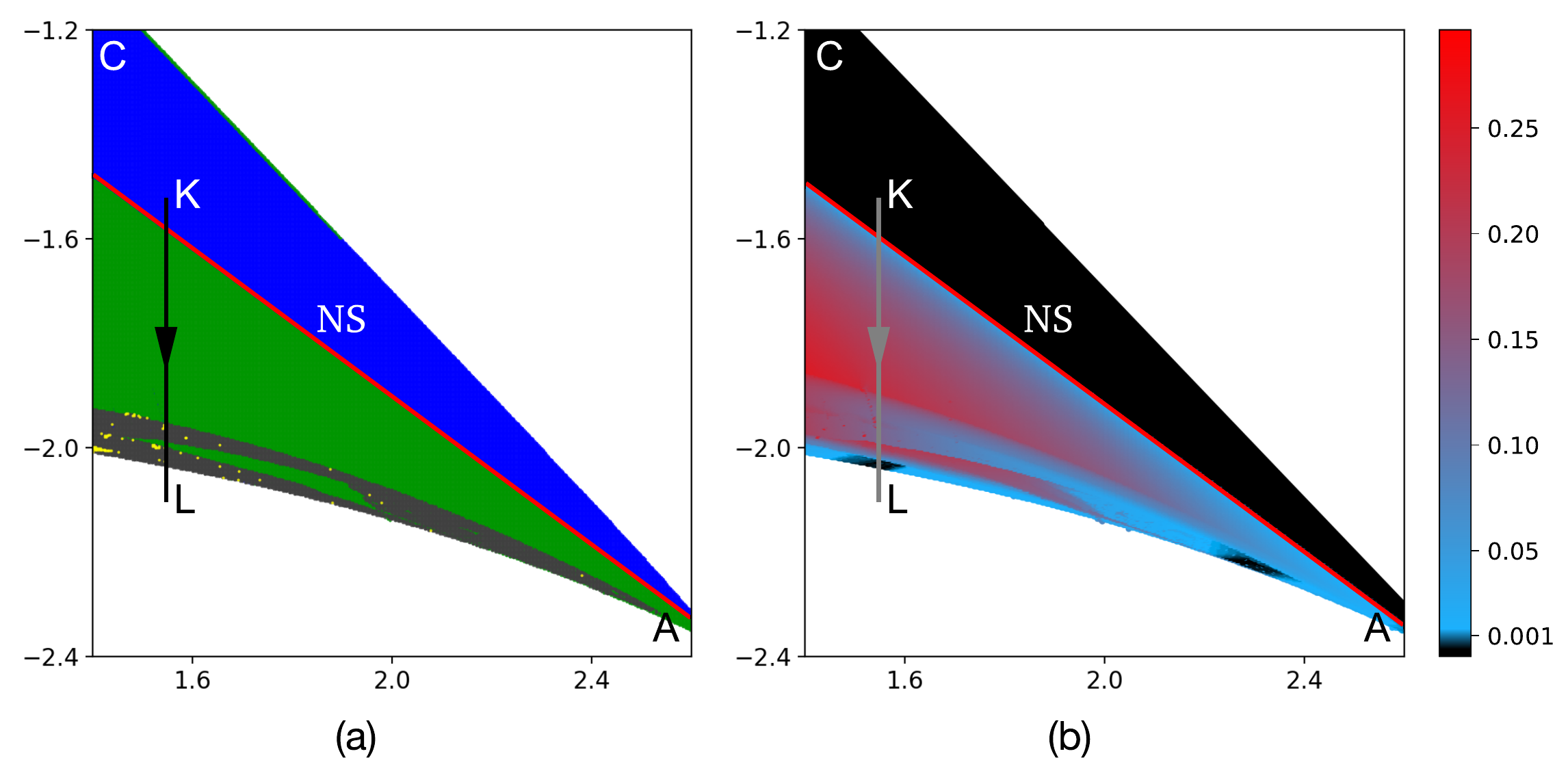}}
\end{minipage}
\caption{\footnotesize ($A,C$)-parameter diagrams for map \eqref{eq:Henon3DMap}, $B=0.7$. (a) Lyapunov diagram, the color scheme is the same as in Fig.~\ref{ris:Fig20}a. (b) Diagram of the distance between the attractor and the fixed point $O_+$, in black colored regions this distance is less than $0.001$. Homoclinic attractors containing the saddle-focus point $O_+$ of type (1,2) appear in the bottom regions, colored in black in panel (b). Almost all chaotic attractors observed here have one positive and one near-zero Lyapunov exponents.}
\label{ris:Fig23}
\end{figure}

Let us consider one more case for the appearance of a flow-like Shilnikov attractor in map \eqref{eq:Henon3DMap}. Now we take $B=0.7$, see the Lyapunov diagram in Fig.~\ref{ris:Fig10}d and its enlarged fragment near the right-bottom corner in Fig.~\ref{ris:Fig23}a. The corresponding distance diagram is presented in Fig.~\ref{ris:Fig23}b. Analyzing these figures, one can see that there are two regions where homoclinic attractors exist. However, non of these regions admit hyperchaotic attractors. Moreover, hyperchaotic attractors are not observable in this case at all, we find only flow-like chaotic attractors.

\begin{figure}[h!]
\begin{minipage}[h]{1\linewidth}
\center{\includegraphics[width=1\linewidth]{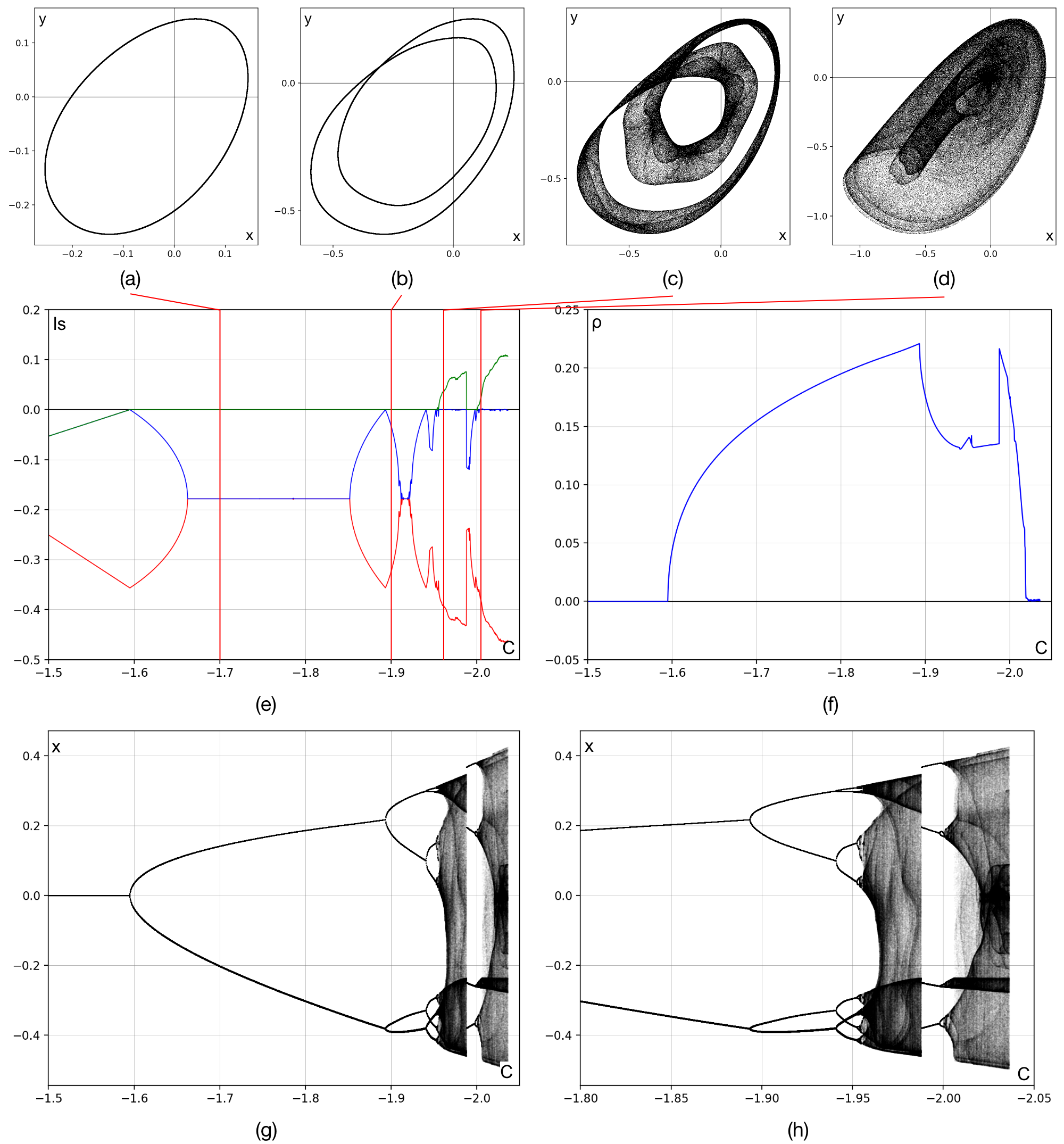}}
\end{minipage}
\caption{\footnotesize Illustration of the attractor evolution along the pathway KL: $(B=0.7, A=1.55)$. (a)--(d) phase portraits of attractors: (a) $C=-1.7$ -- the stable invariant curve $L$;(b) $C=-1.9$ -- doubled invariant curve; (c) $C=-1.96$ -- torus-chaos attractor (LE: $\Lambda_1 = 0.0348$, $\Lambda_2 = 0.0001$, and $\Lambda_3 = -0.392$); (d) $C=-2.03$ -- discrete Shilnikov attractor (LE: $\Lambda_1 = 0.107$, $\Lambda_2 \approx 0$, and $\Lambda_3 = -0.463$). (e) The graph of Lyapunov exponents $\Lambda_1, \Lambda_2$, and $\Lambda_3$ on parameter $C$. (f) the graph of the distance between the attractor and the fixed point $O_+$. (g) Bifurcation tree $x(C)$ computing with help of a Poincar\'e-like map on the plane $(x,z)$ (using the cross-section-like box $|y| < 0.001$) and (f) its enlarged fragment.}
\label{ris:Fig24}
\end{figure}

Let us consider the pathway KL: $(B=0.7, A=1.55)$ and study bifurcations leading to the birth of homoclinic attractors in the left-bottom part of the diagrams presented in Fig.~\ref{ris:Fig23}. Some results of the corresponding one-parameter bifurcation analysis are shown in Fig.~\ref{ris:Fig24}. In this case, the scenario for the appearance of discrete Shilnikov attractors is almost the same as it was in the previous case, cf. Figs.~\ref{ris:Fig24}a--\ref{ris:Fig24}f with Figs.~\ref{ris:Fig21}a--\ref{ris:Fig21}f. The main difference between these two cases is that the stable invariant curve $L$ undergoes here a quite long sequence of the period-doubling bifurcations before destruction of the resulting multiround invariant curve.

Fig.~\ref{ris:Fig24}g shows the corresponding bifurcation tree which we compute using the following scheme. We construct a Poincar\'e-like maps on the plane ($x,z$) introducing a cross-section box $|y| < 0.001$ in the phase space of the map. If a point of the attractor falls into this box we take its $x$-coordinate and plot on the graph $x(C)$. For each value of the parameter $C$ we iterate the map 10000 times and plot the last 3000 points. The enlarged fragment of this bifurcation tree is shown in Fig.~\ref{ris:Fig24}h. It is iteresting to note that the observed bifurcation tree looks like the well-known Feigenbaum tree \cite{Feig79} accompanying the formation of H\'enon-like attractors. As it is known, such attractors often appear after a cascade of period-doubling bifurcations followed by a cascade of heteroclinic band-fusion bifurcations which result in the absorption of periodic saddle orbits emerging after the corresponding period-doubling bifurcations. By the same manner, a sequence of the period-doubling bifurcations with a stable invariant curve leads to the appearance of multicomponent chaotic attractor (on the corresponding two-dimensional Poincar\'e-like map) which, then, transforms to the one-component attractor via a sequence of heteroclinic bifurcations resulting in the absorption of the saddle invariant curves emerging after the corresponding period-doubling bifurcations. The resulting one-component attractor were called in \cite{BSV02, BSV05, BSV10} quasiperiodic H\'enon-like attractor. It has one positive and one near-zero (indistinguishable from zero in numerics) Lyapunov exponents.

Finally, the chaotic attractor absorbs the saddle-focus fixed point $O_+$ of type (1,2), and a discrete Shilnikov flow-like attractor appears, see Fig.~\ref{ris:Fig24}d. The corresponding homoclinic structure for $O_+$ is shown in Fig.~\ref{ris:Fig22}b.

A natural questions arise here.
\begin{itemize}
\item Why does the second Lyapunov exponent is near-zero for chaotic attractors of the three-dimensional map inside large open regions of the parameter space?
\item Why do these attractors look like chaotic attractors of some three-dimensional system of ODEs, if the points of an orbit taken in the attractor are drawn not too densely?
\end{itemize}

\begin{figure}[h]
\center{\includegraphics[width=0.99\linewidth]{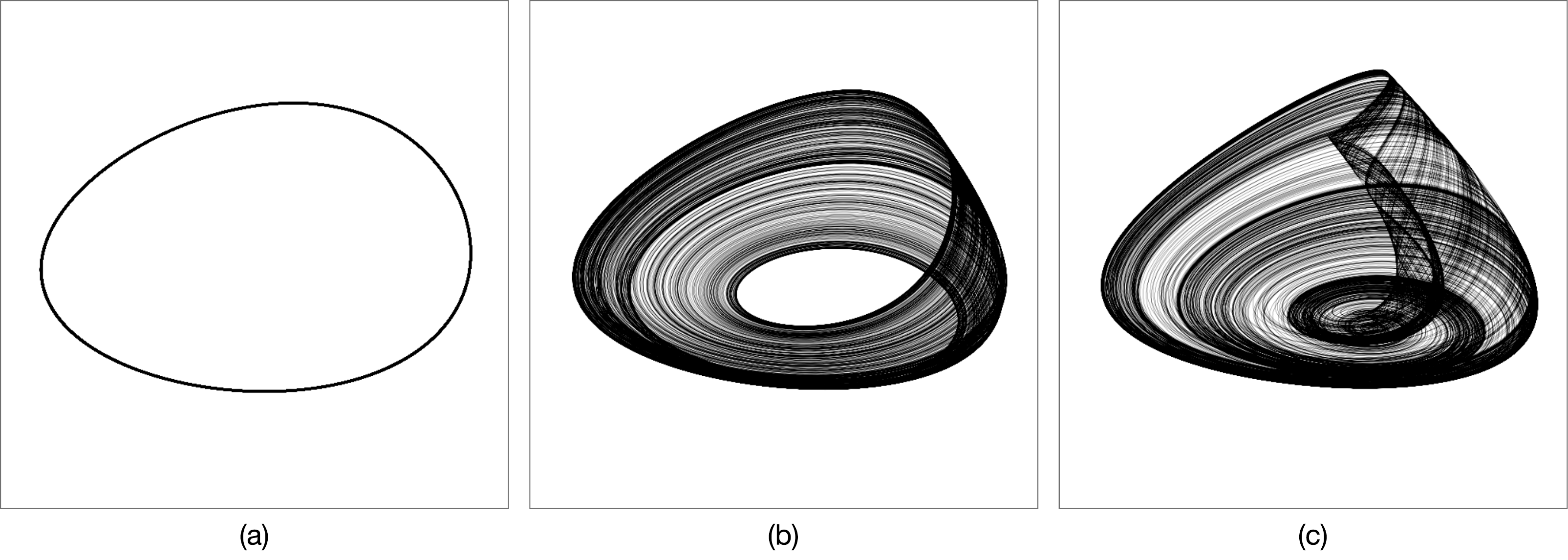}}
\caption{{\footnotesize Main steps in the scenario of the Shilnikov attractor appearance in system \eqref{eq_ACT} for $C=-0.4, B=-1$: (a) $A=-0.7$, stable limit cycle, (b) $A=-0.8$, R\"ossler-like attractor, (c) $A=-0.87$, Shilnikov spiral attractor.}}
\label{Fig25}
\end{figure}

The answer to these questions is quite simple. All flow-like chaotic attractors in map \eqref{eq:Henon3DMap} are observed close enough to the codimension-three bifurcation, when fixed point $O_+$ has the triplet of multipliers $(1,1,1)$: a pair of multipliers $\mu_{1,2}$ are equal to $1$ at the resonance $1:1$ bifurcation, where the curve NS intersects with the curve TR, while the third multiplier $\mu_3 = J$ is equal to one when the map is conservative (volume-preserving). It is not difficult to show that the asymptotic normal form for this codimension-three bifurcation in the case of map \eqref{eq:Henon3DMap} coincides with the well-known Arneodo-Coullet-Spiegel-Tresser (ACST) system
\begin{equation}
\begin{cases}
\dot x = y, \\
\dot y = z, \\
\dot z = \alpha x + \beta y + \gamma z + x^2.
\label{eq_ACT}
\end{cases}
\end{equation}
Bifurcation analysis for this system were performed in detail in Refs.~\cite{ACST1985,ACST1985b}. In particular, it was shown that Shilnikov attractors containing saddle-focus equilibrium of type (1,2) with a homoclinic loop exist in this system. In Refs.~\cite{GGKKB19,bakhanova2022shilnikov} it was described how such attractors are born from the stable equilibrium via chain of local and homoclinic bifurcations. At first, a stable equilibrium $Eq(0,0,0)$ undergoes the supercritical Andronov-Hopf bifurcation after which it becomes a saddle-focus and a stable limit cycle is born in its neighborhood, see Fig.~\ref{Fig25}a. Then, this cycle goes through the cascade of period-doubling bifurcations, see Fig.~\ref{Fig25}b. With further increase in the control parameter, orbits of the attractor come closer and closer to the saddle-focus and, finally, $Eq$ is included to the attractor, as a result of the appearance of Shilnikov homoclinic loop, see Fig.~\ref{Fig25}c. As any chaotic attractor of an ODE system, this attractor has one positive and one zero Lyapunov exponents. This property is inherited by map \eqref{eq:Henon3DMap} in some (sufficiently large) neighborhood of the codimension-three bifurcation when fixed point $O_+$ has the triplet of multipliers $(1,1,1)$.\footnote{Similar phenomenon of the existence of flow-like chaotic attractors was previously discovered and studied in Ref.~\cite{GOST05} for a 3D H\'enon map $\bar x = y, \; \bar y = z, \;  \bar z = M_1 + Bx + M_2y - z^2$ possessing the discrete Lorenz attractor (with zero second Lyapunov exponent) in an open large region of the ($M_1,M_2,B$)-parameter space adjacent to the codimension-three point when a fixed point of the map has multipliers $(-1,-1,1)$. The 3D H\'enon map appears as a normal form for the first-return map near a quadratic homoclinic tangency to a saddle-focus fixed point with the Jacobian closed to 1~\cite{GMO06}. Also note that this map is the inverse of the 3D Mir\'a map \eqref{eq:MiraMap}. Therefore, in map \eqref{eq:MiraMap} there exists a discrete Lorenz repeller (with zero second Lyapunov exponent) for $B > 1$.}

If to compare phase portraits of the attractor found in map \eqref{eq:Henon3DMap} at $B=0.7$ (Fig.~\ref{ris:Fig22}b) and the attractor presented in Fig.~\ref{Fig25}c, it is easy to conclude that these attractors are similar. An orbit taken in the first attractor looks like a discretization of a continuous orbit taken in the second attractor.

\section{Discussion}

We have proposed bifurcation scenarios leading from a stable fixed point to a hyperchaotic attractor in one-parameter families of three-dimensional maps and have applied them for studying the formation of hyperchaotic attractors in the homoclinic three-dimensional Mir\'a map~\eqref{eq:MiraMap}. These scenarios consist of two parts. The first, phenomenological, part describes a few main bifurcations resulting to the appearance of a homoclinic attractor containing saddle fixed point of type (1,2). The second, empirical, part describes accompanying bifurcations after which the majority of orbits inside an attractor get two-dimensional unstable manifolds.
We have shown that for map~\eqref{eq:MiraMap} phenomenological part is the extension of the well-known Shilnikov scenario \cite{Sh86} to the case of three-dimensional maps \cite{GGS12,GGKT14,GG16}. In the framework of this scenario, a stable fixed point $O_+$ undergoes the supercritical Neimark-Sacker and a stable invariant curve $L$ appears in its neighborhood. Then, this curve breaks down giving a torus-chaos attractor with the isolated saddle-focus point $O_+$ of type (1,2). The final part of the scenario is the inclusion (absorption) of this point by the attractor. The resulting discrete Shilnikov attractor contains this point.

\begin{figure}[h!]
\begin{minipage}[h]{1\linewidth}
\center{\includegraphics[width=0.8\linewidth]{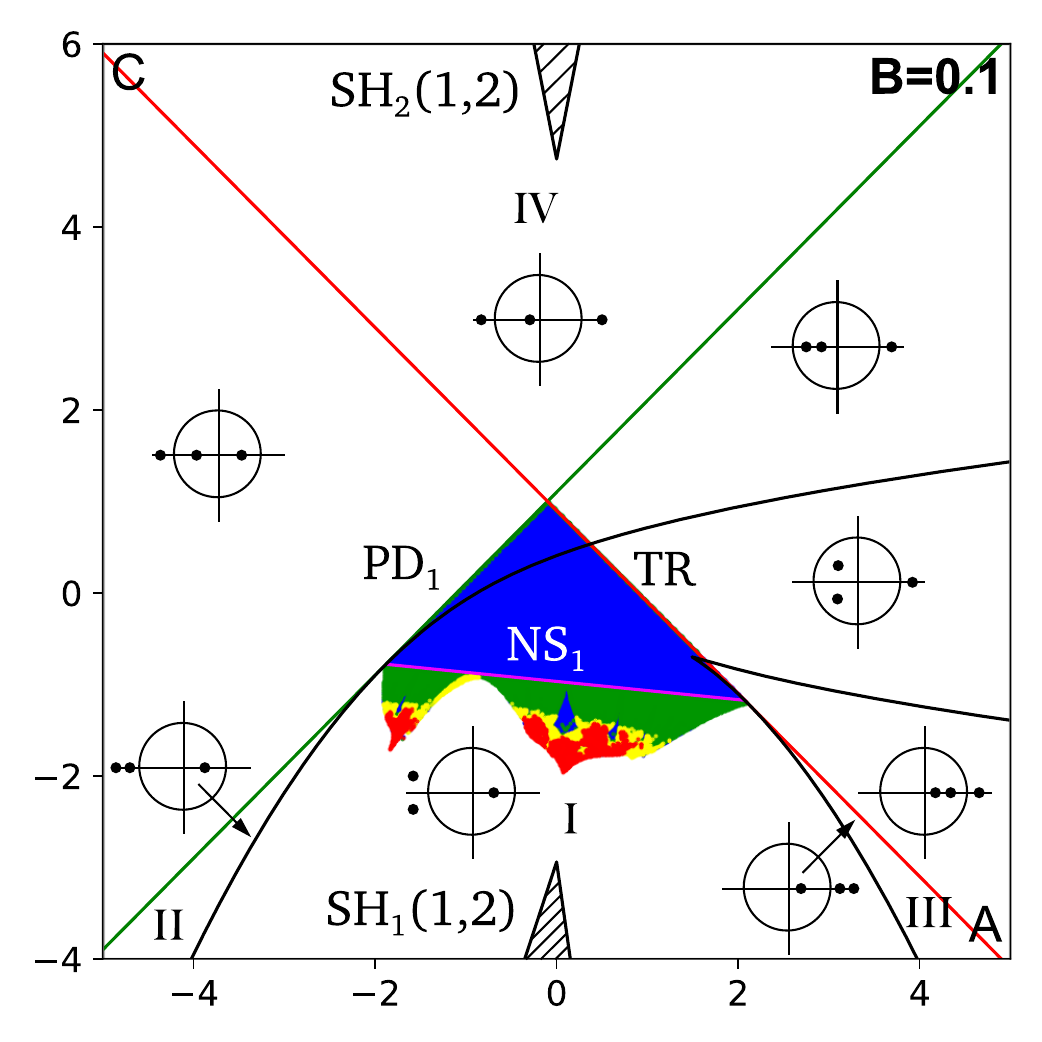}}
\end{minipage}
\caption{\footnotesize Saddle chart superimposed with the Lyapunov diagram for map \eqref{eq:Henon3DMap}, $B=0.1$. In regions $SH_{1,2}(1,2)$ the nonwandering set of the map is nontrivial hyperbolic of type (1,2). Hyperchaotic homoclinic attractors of different types are possible in regions I--IV.}
\label{ris:Fig26}
\end{figure}

The empirical part of the scenario describes mechanisms of the destruction of the curve $L$. Usually, before the destruction, this curve becomes resonant. We have shown that depending on values of parameters the corresponding resonant orbits can give rise to hyperchaotic periodic orbits of type (1,2) in two different ways:
\begin{enumerate}[label=(\roman*)]
\item the stable resonant orbit undergoes a cascade of period-doubling bifurcations. (This cascade can be interrupted by the supercritical Neimark-Sacker bifurcation. In this case, see (ii)). In their turn, the resulting periodic saddle orbits of type (2,1), as well as the resonant saddle orbits of type (2,1) undergo a cascade of period-doubling bifurcations transforming the corresponding periodic orbits to saddles of type (1,2).
\item the stable resonant orbit undergoes the supercritical Neimark-Sacker bifurcation transforming this orbit to a saddle-focus of type (1,2) and a stable multicomponent invariant curve appears. Then, this invariant curve (as the curve $L$) becomes resonant, (i) or (ii) is implemented and so on. The saddle resonant orbit, as well as in the first way, undergoes a cascade of period-doubling bifurcations.
\end{enumerate}
In both cases, an attractor absorbs the sets of periodic saddle orbits of type (1,2) and becomes hyperchaotic.

In the paper, we have also proposed a new phenomenological scenario leading to the creation of the so-called hyperchaotic H\'enon-like attractor containing a saddle fixed point with a pair of negative unstable multipliers, However, we have not observed its implementation in map \eqref{eq:Henon3DMap}.

Finally, we would like to note one important property of such maps. One of its fixed points is always at the origin and its eigenvalues depend only on the parameters $A,B$, and $C$ \cite{GG16}. In Figure~\ref{ris:Fig26}, we show an extended bifurcation diagram for this fixed point (when $B=0.1$) above the Lyapunov diagram. Such extended diagrams were called \textit{saddle charts} in \cite{GG16}.

In this diagram one can see four regions with possible hyperchaotic homoclinic attractors: the region I -- with a discrete Shilnikov attractor, the region II -- with a hyperchaotic H\'enon-like attractor, the region III -- with a hyperchaotic attractor containing a fixed point with a pair of positive unstable multipliers and the region IV -- with a hyperchaotic attractor containing a fixed point with a pair of real unstable multipliers with different signs. Varying nonlinear terms in the map we can expect new types of hyperchaotic homoclinic attractors. This problem looks very promising, especially for the region IV, inside which, as well as inside the region I, by Gonchenko and Li theorem \cite{GonLi10}, the nonwandering set of the map is a Smale horseshoe of type (1,2). Thus, one can expect the existence of hyperchaotic attractors in IV for quite general families of three-dimensional H\'enon maps.

\vspace{0.5cm}
\textbf{Acknowledgements}

The authors thank Sergey Gonchenko, Igor Sataev, and Dmitry Turaev for fruitful discussions. The authors also express their gratitude to the referee for valuable comments, which significantly helped to improve this paper.

This work was carried out with the support of the RSF grant No.~19-71-10048 (Sections~\ref{sec:SmallB} and \ref{sec:LargeB}) and the Laboratory of Dynamical Systems and Applications NRU HSE, grant for supporting basic scientific research in inter-campus groups (Sections~\ref{sec:Scenarios}, \ref{sec:MiraMapBiff}, and \ref{sec:VarShil}). A.~Kazakov and E.~Karatetskaia also acknowledge the support of Theoretical Physics and Mathematics Advancement Foundation BASIS for financial support of scientific investigations.	

\newpage

\printbibliography

\end{document}